\RequirePackage{fix-cm}
\documentclass[leqno, fleqn, centertags, 12pt]{article}

\usepackage{latexsym}
\usepackage{amsmath}
\usepackage{amscd}
\usepackage{amssymb}
\usepackage{verbatim}

\usepackage[T1]{fontenc}
\usepackage{latexsym}
\usepackage{manfnt}
\usepackage[sans]{dsfont}
\usepackage{palatino}
\usepackage[sc, osf]{mathpazo} 
\usepackage{eulervm}
\usepackage{stmaryrd}

\usepackage{fullpage}

\usepackage{graphicx}

\DeclareMathOperator{\Aut}{{Aut}\hspace{0.5pt}}
\DeclareMathOperator{\ccc}{{\mathcal C}\hspace{0.5pt}} 
\DeclareMathOperator{\genus}{{\hspace*{-0.1em}{\mathfrak g}\hspace{-0.1em}}}
\DeclareMathOperator{\Inn}{Inn\hspace{0.5pt}}
\DeclareMathOperator{\inter}{int\hspace{0.5pt}}
\DeclareMathOperator{\Mod}{{Mod}\hspace{0.5pt}}
\DeclareMathOperator{\Out}{Out\hspace{0.5pt}}
\DeclareMathOperator{\rank}{{rank}\hspace{0.5pt}}
\DeclareMathOperator{\teich}{T\hspace{0.5pt}}
\DeclareMathOperator{\tei}{T}

\newcommand{\mm}{\Mod_{\hff\genus}}
\newcommand{\xx}{X_{\hff\genus}}
\newcommand{\ccg}{\ccc(\fff\xx\fff)}
\newcommand{\cht}{\mathcal{HT}\hspace{0.5pt}(X_{\hff\genus})}
\newcommand{\chtone}{\mathcal{HT}_1\hspace{0.5pt}(X_{\hff\genus})}
\newcommand{\vcd}{\mathrm{vcd}\hspace*{0.1em}}
\newcommand{\ccd}{\mathrm{cd}\hspace*{0.1em}}
\newcommand{\ttg}{\teich_{\hff\genus}}
\newcommand{\ttb}{\overline{\tei}\hspace{0.5pt}_{\hff\genus}}

\newcommand{\xxb}{\overline{X}}

\newcommand{\rr}{\mathbold{R}}

\newcommand{\zz}{\mathbold{Z}}

\newcommand{\card}{\mathop{\mathrm{card}}\nolimits}

\newcommand{\tto}{\,{\to}\,}

\newcommand{\langles}{\mathbin{\langle}}  
\newcommand{\rangles}{\mathbin{\rangle}}  
\newcommand{\vv}[1]{\langles #1 \rangles} 
\newcommand{\llvv}[1]{\langles\hspace*{0.05em}{#1}\hspace*{0.07em}\rangles}
\newcommand{\vvv}[1]{\hspace*{-0.2em}\langles\hspace*{-0.15em}{#1}\hspace*{0.05em}\rangles}

\newcommand{\geqs}{\hskip.05em\,{\geq}\,}      
\newcommand{\leqs}{\,{\leq}\,\hskip.05em}        
\newcommand{\eeq}{{\hskip.05em\,{=}\,}}         
\newcommand{\nneq}{{\hskip.05em\,{\neq}\,}}      

\newcommand{\gres}{\hskip.05em\,{>}\,}          
\newcommand{\less}{\,{<}\,\hskip.05em}

\newcommand{\setmiss}{{\setminus}\hspace{0.05em}}  
  
\newcommand{\aasetmiss}{\,\setmiss\,}

\newcommand{\ssub}{\,{\subset}\,}

\newcommand{\elem}{\,{\in}\,}

\def\hss{\hskip.025em\ }
\def\sss{\hskip.05em\ }
\def\dss{\hskip.1em\ }
\def\trs{\hskip.15em\ }
\def\qss{\hskip.02em\ }

\def\hff{{\hskip.025em}}
\def\fff{{\hskip.05em}}
\def\dff{{\hskip.1em}}
\def\trf{{\hskip.15em}}
\def\qff{{\hskip.2em}}

\def\endss{\hspace*{0.05em}}

\def\ffdot{\hspace*{-0.1em}.\hspace*{0.1em}\ }
\def\dfdot{\hspace*{-0.2em}.\hspace*{0.2em}\ }

\def\ffcom{\hspace*{-0.1em},\hspace*{0.1em}\ }
\def\dfcom{\hspace*{-0.2em},\hspace*{0.2em}\ }

\newcommand{\nsp}{\hspace*{-0.1em}}
\newcommand{\dnsp}{\hspace*{-0.2em}}

\parskip=\the\bigskipamount

\renewcommand{\refname}{\textnormal{Original references}}

\vspace{-1.5ex}

\makeatletter
\renewcommand{\@makefntext}[1]{\vspace*{0.5ex}\parindent=0em\noindent
\hspace*{-0.4em}
\hbox to 0.4em{\hss\@makefnmark}\hspace*{0.4em}{#1}
}
\makeatother

\newcounter{mysectionnumber}
\newcounter{myparnum}
\newcounter{mylemmanum}[myparnum]
\newcommand{\mysection}[2]{\setcounter{footnote}{0}
\refstepcounter{mysectionnumber}
\section*{ \textnormal{{\themysectionnumber.} {#1}}}\label{#2}\vspace{\bigskipamount}}

\newcommand{\mynonumbersection}[2]{
\vspace{-2.0ex}
\section*{ \textnormal{{#1}}}\label{#2}\vspace{\bigskipamount}\vspace{-2.0ex}}

\newcommand{\mypar}[2]{\refstepcounter{myparnum}
{\vspace{-\bigskipamount} \paragraph{\textit{{
\themyparnum. #1}\label{#2}}} \hspace{-0.5em}}}
\renewcommand{\themyparnum}{\arabic{myparnum}}

\newcounter{mysubsubnumber}[myparnum]
\newcounter{mysubparnumber}[myparnum]
\newcommand{\myitpar}[1]{\vspace{-\bigskipamount}\paragraph{\textit{#1}}\hspace{-0.7em}}

\newcommand{\proof}{\vspace{-\bigskipamount}{\paragraph{{\emph{Proof}.\hspace*{0.2em}}}\hspace{-0.7em}}}
\newcommand{\eproof}{ $\blacksquare$}
\newcommand{\esubproof}{ $\square$}

\newcommand{\minus}{\hspace*{0.15em}\mbox{\rule[0.4ex]{0.4em}{0.4pt}}\hspace*{0.15em}}
\newcommand{\mminus}{\hspace*{0.15em}\mbox{\rule[0.5ex]{0.4em}{0.4pt}}\hspace*{0.15em}}
\newcommand{\plus}{\hspace*{0.15em}\mbox{\rule[0.4ex]{0.5em}{0.4pt}\hspace*{-0.27em}\rule[-0.15ex]{0.4pt}{1.2ex}\hspace*{0.27em}}\hspace*{0.15em}}
\newcommand{\pplus}{\hspace*{0.15em}\mbox{\rule[0.5ex]{0.5em}{0.4pt}\hspace*{-0.27em}\rule[-0.1ex]{0.4pt}{1.2ex}\hspace*{0.27em}}\hspace*{0.15em}}

\begin{document}

\title{The\hspace*{0.1em} virtual\hspace*{0.1em} cohomology\hspace*{0.1em} dimension\hspace*{0.1em}\\  
of\hspace*{0.1em} Teichm{\"u}ller\hspace*{0.1em} modular\hspace*{0.1em} groups\hspace*{0.1em}:\hspace*{0.3em}\\
the\hspace*{0.1em} first\hspace*{0.1em} results\hspace*{0.1em} and\hspace*{0.1em}
a\hspace*{0.1em} road\hspace*{0.1em} not\hspace*{0.1em} taken
\vspace{2.5ex}
\begin{flushright}
\large
\textit{To Leningrad Branch of Steklov Mathematical Institute, \\ 
with gratitude for giving me freedom to \\ pursue my interests in 1976--1998.}
\end{flushright}
\vspace{0.3ex}
\vspace*{1.5ex}} 
\date{}
\author{\textnormal{Nikolai\hspace*{0.1em} V.\hspace*{0.1em} Ivanov{}\vspace*{6.3ex}}}

\flushbottom

\maketitle

\footnotetext{\hspace*{-0.5em}\copyright\ Nikolai V. Ivanov, 2015.\trs 
The results of Sections \ref{introduction}, \ref{boundary}, and \ref{ccht}
were obtained in 1983 at the Leningrad Branch of Steklov Mathematical Institute.
The preparation of this paper was not supported by any governmental or non-governmental agency, foundation, or institution.}

\mynonumbersection{Preface}{preface}

{\small 
Sections \ref{introduction}, \ref{boundary}, and \ref{ccht} of this paper
are based on my paper\sss \cite{iv-vcd}.\sss
Section \ref{motivation} is devoted to the motivation behind the results of \cite{iv-vcd}.\sss 
Section \ref{context} is devoted to the context of my further results about the virtual cohomological dimension
of Teichm{\"u}ller modular groups.
The last Section \ref{myremin} is devoted to the mathematical and non-mathematical circumstances 
which shaped the paper\sss \cite{iv-vcd} and some further developments. 

As I only recently realized,\sss \cite{iv-vcd}\sss contains the nucleus of some techniques 
for working with complexes of curves and other similar complexes 
used many times by myself and then by other mathematicians.
On the other hand, one of the key ideas of\sss \cite{iv-vcd},\sss namely, 
the idea of using the Hatcher--Thurston cell complex \cite{ht},
was abandoned by me already in 1983, 
and was not taken up by other mathematicians.
It seems that it still holds some promise.
This is a \emph{road not taken}, alluded to in the title.
The implied reference to a famous poem by Robert Frost is indeed relevant, 
especially if the poem is understood not in the clich\'{e}d way.

The list of references consists of two parts. 
The first part reproduces the list of references of\sss \cite{iv-vcd}.\sss 
The second part consists of additional references.
The papers from the first list are referred to by numbers, and from the second one by letters followed by numbers.
So, \cite{bs} refers to the first list, and \cite{iv-vcd} to the second.

The exposition and English in Sections \ref{introduction}, \ref{boundary}, and \ref{ccht} are, 
I hope, substantially better than in\sss \cite{iv-vcd}.\sss 
At the same time, these sections closely follow\sss \cite{iv-vcd}\sss with one exception.
Namely, the original text of\sss \cite{iv-vcd}\sss contained a gap.
A densely written correction was added to\sss \cite{iv-vcd}\sss at the last moment as an additional page.
In the present paper this correction is incorporated into the proof of Lemma \ref{loops} (see the subsection \ref{proof-loops}, Claims 1 and 2). 
Finally, \LaTeX\ leads to much better output than a typewriter combined with writing in formulas by hand.

}

\mysection{Introduction}{introduction}

Let $X_g$ be a closed orientable surface of genus ${\genus}\geqs 2$\dfdot
The \emph{Teichm{\"u}ller modular group} $\mm$ of genus $\genus$ is defined as the group of isotopy classes of diffeomorphisms
$\xx\tto\xx$,\sss i.e. $\mm\dnsp\eeq\fff\pi_0(\hff \mbox{Diff}\dff(\xx))$\dfdot
This group may be also defined as the group of homotopy classes of homotopy equivalences $\xx\tto\xx$\ffdot
These two definitions are equivalent by a classical result of Baer--Nielsen.
See \cite{zi} for a modern exposition close in spirit to the original papers of J. Nielsen and R. Baer.
It is well known that $\xx\fff$ is a $K(\pi,1)$\dnsp-space, where $\pi\eeq\pi_1(\xx)$\dfdot
This allows to determine the group of homotopy classes of homotopy equivalences $\xx\tto\xx$ in terms of the fundamental group $\pi_1(\xx\fff)$ alone 
and conclude that $\mm$ is isomorphic to the outer automorphisms group
\[
\Out(\pi_1(\xx\fff)) = \Aut(\pi_1(\xx\fff))/\Inn (\pi_1(\xx\fff))\endss.
\] 
Here $\Aut(\pi)$ denotes the automorphism group of a group $\pi$ and $\Inn(\pi)$ denotes the subgroup of inner automorphisms of $\pi$\dfdot
The groups $\mm$ are also known as the \emph{surface mapping class groups}.

The present paper is devoted to what was the first step toward the computation of the virtual cohomology dimension of the groups $\mm$.
Its main result (Theorem \ref{vcd} below) provided the first non-trivial estimate of the virtual cohomology dimension of $\mm$.

The ordinary cohomology dimension of $\mm$ is infinite because $\mm$ contains non-trivial elements of finite order.
However $\mm$ is virtually torsion free, i.e. $\mm$ contains torsion free subgroups of finite index
(this result is due to J.-P. Serre; see \cite{gr}\fff).
If a group $\Gamma$ is virtually torsion free, then all torsion free subgroups of $\Gamma$ of finite index
have the same cohomology dimension, which is called the \emph{virtual cohomology dimension} of $\Gamma$ and is denoted by $\vcd(\Gamma)$\dfdot
It is well known that $\vcd(\mm)$ is finite and moreover that $\vcd(\mm)\eeq 1$ for $\genus\eeq 1$ and
 $3\genus{\minus}3\leqs\vcd(\mm)\leqs 6\genus {\minus}7$ for $\genus\geqs 2$\dfdot

\begin{figure}[t]
\includegraphics[width=0.96\textwidth]{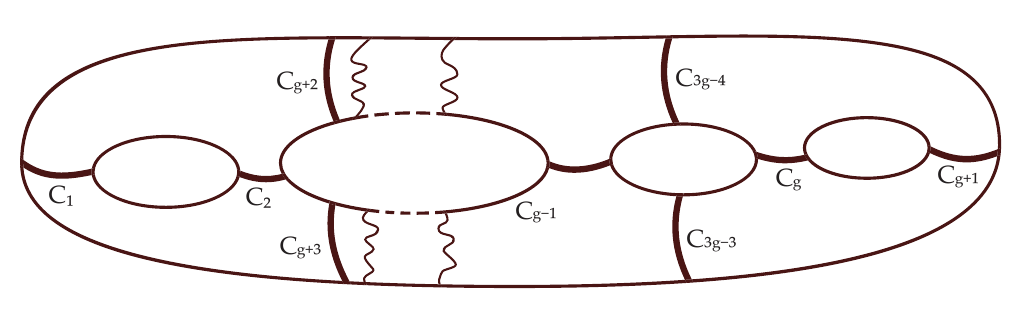}
\caption{A maximal collection of disjoint pair-wise non-isotopic circles.}
\label{3g-3}
\end{figure}

Let me recall the proofs of the last inequalities (assuming that $\genus\geqs 2$\nsp). 
In order to prove the first one, recall that $\mm$ contains
a free abelian subgroup of $\rank = 3\genus{\minus}3$\dfdot
For example, the subgroup generated by Dehn twists along the curves $C_1\dff,\dff\ldots\dff,\dff C_{3\genus{\minus}3}$
on Fig. \ref{3g-3} is a free abelian subgroup of $\rank = 3\genus{\minus}3$\dfdot
Since the virtual cohomology dimension cannot be increased by passing to a subgroup,
and since $\vcd(\zz^n)\eeq n$\dfcom
we see that $3\genus{\minus}3\leqs\vcd(\mm)$\dfdot

\renewcommand{\topfraction}{1.0}
\renewcommand{\textfraction}{0.0}

The proof of the inequality $\vcd(\mm)\leqs 6\genus {\minus}7$ is much more deep and 
is based on theories of Riemann surfaces and of Teichm{\"u}ller spaces.
Recall that $\mm$ naturally acts on the Teichm{\"u}ller space $\ttg$ of marked Riemann surfaces of genus $\genus$\dfdot
The action of $\mm$ on $\ttg$ is a properly discontinuous, and the quotient space $\ttg/\nsp\mm$ is the \emph{moduli space} of Riemann surfaces
of genus $\genus$ (\hff this is the source of the term \emph{Teichm{\"u}ller modular group}\fff).
Moreover, any torsion free subgroup $\Gamma$ of $\mm$ acts on $\ttg$ freely{}.\qss
Since $\ttg$ is homeomorphic to $\rr^{6\genus{\minus}6}$\dfcom
for such a subgroup $\Gamma$ the quotient space $\ttg/\Gamma$ is a $K(\Gamma,1)$\dnsp-space (and, in addition, is a manifold).
This implies that the cohomology dimension of $\Gamma$ is $\leqs\dim \ttg/\Gamma\eeq 6\genus{\minus}6$\dfcom
and hence $\vcd(\mm)\leqs 6\genus{\minus}6$\dfdot
In order to prove that, moreover, $\vcd(\mm)\leqs 6\genus{\minus}7$\dfcom
recall that $\ttg/\nsp\mm$ is non-compact, and hence $\ttg/\Gamma$ is also non-compact.
Since the $n$\dnsp-th cohomology groups of $n$\dnsp-dimensional non-compact manifolds with any coefficients,
including twisted ones, are equal to $0$, this implies that the $6\genus{\mminus}6$\dnsp-th 
cohomology group of any such subgroup $\Gamma$ is equal to $0$\dfcom
and hence the cohomology dimension of $\Gamma$ is $\less 6\genus{\minus}6$\dfdot
It follows that $\vcd(\mm)\less 6\genus{\minus}6$\dfcom
i.e. $\vcd(\mm)\leqs 6\genus{\minus}7$\dfdot 

The main result of the present paper is the following strengthening of the inequality $\vcd(\mm)\leqs 6\genus{\minus}7$\dfdot

\mypar{Theorem.}{vcd} $\vcd(\mm)\leqs 6\genus{\minus}9$ 
\emph{for $\genus\geqs 2$ and $\vcd(\Mod_2)\eeq 3$\dfdot
In addition, $\Mod_2$ is virtually a duality group of dimension $3$\dfdot}

The proof of this theorem is based on the properties of a boundary of Teichm{\"u}ller space introduced by W. Harvey \cite{ha1, ha2}.
The key property of Harvey's boundary of $\ttg$ is the fact that it is homotopy equivalent to (the geometric realization of) a simplicial complex 
$\ccg$\dfdot
The complex $\ccg$ was also introduced by W. Harvey and is known as the \emph{complex of curves} of $\xx$\nsp.\dss
We will recall the definition of complexes of curves in Section 3.
Using the results of W. Harvey and the theory of cohomology of groups, 
especially the theory of groups with duality developed by R. Bieri and B. Eckman
Theorem \ref{vcd} can be reduced to the following theorem (see Section \ref{boundary}).

\mypar{Theorem.}{ccsimply} \emph{The complex of curves $\ccg$ of\dss $\xx$ is simply-connected for\sss $\genus\geqs 2$\dfdot}

Using the same arguments one can deduce from Theorem \ref{ccsimply} 
that $\Mod_2$ is virtually a duality group in the Bieri-Eckmann sense \cite{bie},
i.e. that $\Mod_2$ contains a subgroup of finite index which is a duality group in the Bieri-Eckmann sense.
Theorem \ref{ccsimply} is deduced from the simply-connectedness of a cell complex introduced by A. Hatcher and W. Thurston
\cite{ht} (see Section \ref{ccht}).
The simply-connectedness of the Hatcher-Thurston complex is one of the main results of their paper \cite{ht}.

The rest of the paper is arranged as follows. 
In Section \ref{boundary} we review the basic properties of Harvey boundary of Teichm{\"u}ller space,
and then deduce Theorem \ref{vcd} from Theorem \ref{ccsimply}.
In Section \ref{ccht} we start with defining complexes of curves and Hatcher-Thurston complexes,
and then deduce Theorem \ref{ccsimply} from results of A. Hatcher and W. Thurston \cite{ht}.

In Section \ref{motivation} we explain the ideas from the theory of arithmetic groups which served as a motivation for
the approach to the virtual cohomology dimension of $\mm$ outlined above, and for the further work in this direction.
In Section \ref{context} we outline a broad context in which Theorem \ref{ccsimply}
and then stronger results about the connectivity of $\ccg$ were discovered.
Section \ref{myremin} is the last one and is devoted to some personal reminiscences related to these stronger results.
It has grown out of a short summary written by me in Summer of 2007 as a step toward writing 
the expository part of the paper \cite{ivj} by Lizhen Ji and myself.

\mysection{Motivation from the theory of arithmetic groups}{motivation}

\vspace*{\bigskipamount}
\myitpar{The Borel-Serre theory.} Around 1970 A. Borel and J.-P. Serre 
studied cohomology of arithmetic \cite{bs-1970} and $S$\dnsp-arithmetic \cite{bs-1971} groups.
In particular, Borel and Serre computed the virtual cohomology dimension of such groups.
The details were published in \cite{bs} and \cite{bs-1976} respectively.

In outline, Borel and Serre approach is as follows.
Let $\Gamma$ be an arithmetic group.
There is a natural contractible smooth manifold $X$ on which $\Gamma$ acts.
Moreover, $\Gamma$ acts on $X$ properly discontinuously,
and a subgroup of finite index in $\Gamma$ acts on $X$ freely.
For the purposes of computing or estimating $\vcd(\Gamma)$\dfcom
we can replace $\Gamma$ by such subgroup, if necessary, and assume that 
$\Gamma$ itself acts on $X$ freely.
Then the quotient $X/\Gamma$ is a $K(\Gamma,1)$\dnsp-space and one may hope to use it
for understanding the cohomological properties of $\Gamma$\dfdot
Unfortunately, $X/\Gamma$ is usually non-compact.

The first step of the Borel-Serre approach \cite{bs-1970}, \cite{bs} is a construction of a natural compactification of $X/\Gamma$\dfdot
This compactification has the form $\xxb/\Gamma$\dfcom where $\xxb$ is a \emph{smooth manifold with corners} independent of $\Gamma$
and having $X$ as its interior.
As a topological space, a \emph{smooth manifold with corners} is a topological manifold with boundary.
It has also a canonical structure similar to that of smooth manifold with boundary:
while the smooth manifolds with boundary are modeled on products $\rr^n{\times}\rr_{\geq 0}$
(where $n{\plus}1$ is equal to the dimension), 
the smooth manifolds with corners are modeled on products $\rr^n{\times}\rr_{\geq 0}^m$
(where $n{\pplus}m$ is the dimension).

The existence of a structure of smooth manifold with corners on $\xxb$ together with the compactness
of $\xxb$ implies that $\xxb/\Gamma$ admits a \emph{finite} triangulation.
In particular, $\xxb/\Gamma$ is homotopy equivalent to a \emph{finite} $CW$-complex.
This implies that the virtual cohomological dimension  $\vcd(\Gamma)$ \emph{is finite}.
In fact, this implies a much stronger finiteness property of $\Gamma$\dfdot
Namely, $\Gamma$ is a \emph{group of type} $(FL)$\dfcom
{\qff}i.e. there exists a resolution of the trivial $\Gamma$-module $\zz$ 
by \emph{finitely generated free modules} and having\qff\qss \emph{finite length}.

The second step of the Borel-Serre method is an identification of the homotopy type of the boundary $\partial X$\dfdot
Borel and Serre proved that $\partial X$ is homotopy equivalent to the (geometric realization of the) 
Tits building associated with $X$ (or, one may say, with $\Gamma$\nsp).
By a theorem of L. Solomon and J. Tits (see \cite{solomon}, \cite{garland}),
the Tits building is homotopy equivalent to a wedge of spheres.
Moreover, all these spheres have the same dimension,
equal to $r{\minus}1$\dfcom where $r$ is the so-called \emph{rank} of $X$ (or of $\Gamma$\nsp).

The last step in the computation of $\vcd(\Gamma)$ by Borel-Serre method is 
an application of a version of the Poincar\'{e}-Lefschetz duality 
(namely, of the version allowing arbitrary twisted coefficients).
This step uses the fact that $\Gamma$ is a group of type $(FL)$ 
implied by the existence of a structure of a smooth manifold with corners on $\xxb$\dfdot
In fact, it would be sufficient to know that $\Gamma$ is a group of type $(FP)$\dfcom
{\qff}i.e. there exists a resolution of the trivial $\Gamma$-module $\zz$ 
by \emph{finitely generated projective modules} and having\qff\qss \emph{finite length}. 

If $\Gamma$ is only a $S$\dnsp-arithmetic group, there is still a natural
contractible smooth manifold $X$ on which $\Gamma$ acts.
But in this case $\Gamma$ does not act on $X$ properly discontinuously.
In order to overcome this difficulty Borel and Serre \cite{bs-1971}, \cite{bs-1976} multiplied $X$ 
by another topological space $Y$ with a canonical action of $\Gamma$\dfdot
The space $Y$ is not a smooth or topological manifold.
In fact, its topology is closely related to the topology of non-archimedean local fields.
This is the source of the main difficulties in the case of $S$\dnsp-arithmetic groups compared to the arithmetic ones.
These difficulties are technically irrelevant for Teichm\"{u}ller modular groups.
But the fact the original Borel-Serre theory can be applied in a situation
different from the original one was encouraging.

\myitpar{The Bieri-Eckmann theory.} While the Borel-Serre theory served as the motivation,
on the technical level it is easier to use more general results of Bieri-Eckmann \cite{bie}. 

In fact, the last step of the Borel-Serre computation of $\vcd(\Gamma)$ works in a very general situation.
The corresponding general theory is due to R. Bieri and B. Eckmann \cite{bie}, who developed it independently of Borel-Serre.
Bieri and Eckmann \cite{bie} presented a polished theory ready for applications.
In the detailed publication \cite{bs-1976} of their results about $S$\dnsp-arithmetic groups  
Borel and Serre used \cite{bie} when convenient. 
The following easy corollary of Theorem 6.2 of Bieri-Eckman
summarizes the results needed.

\mypar{Theorem.}{bie-duality} \emph{Suppose that a discrete group $\Gamma$ 
acts freely on a topological manifold $\xxb$ of dimension $n$ with boundary $\partial\xxb$.
Suppose that $\xxb/\Gamma$ is homotopy equivalent to a finite $CW$-complex.
If for some natural number $d$ the reduced integral homology groups $H_i(\partial\xxb)$ are equal to $0$ for $i\neq d$ 
and the group $H_d(\partial\xxb)$ is torsion-free, then $\vcd(\Gamma)=n-1-d$\nsp.\dss
If the manifold $\xxb/\Gamma$ is orientable, then $\Gamma$ is a duality group in the sense of\dss \textup{\cite{bie}}\dss
and $\ccd(\Gamma)=\vcd(\Gamma)=n-1-d$\nsp.\dss}\vspace*{\medskipamount}

In Borel-Serre theory $\partial\xxb$ is homotopy equivalent to a bouquet of spheres of the same dimension by the Solomon-Tits theorem.
It follows that $H_i(\partial\xxb)=0$ if $i\neq d$\dfcom
where $d$ is the dimension of these spheres,
and that $H_d(\partial\xxb)$ is a free abelian group. 
In particular, it is torsion-free, and hence Theorem \ref{bie-duality} applies.
The next theorem is not proved by Bieri-Eckmann \cite{bie},
but is very close to Theorem \ref{bie-duality}.

\mypar{Theorem.}{vcd-conn-general} \emph{In the framework of Theorem \ref{bie-duality},
if $c$ is a natural number such that 
the reduced homology groups $H_i(\partial\xxb)=0$ for $i\leq c-1$, 
then $\vcd(\Gamma)\leq n-1-c$.}\vspace*{\medskipamount}

Since in the present paper we prove only upper estimates of the virtual cohomology dimension,
Theorem \ref{vcd-conn-general} is better suited for our goals.

\mysection{The Harvey boundary of Teichm{\"u}ller space}{boundary}

An analogue for Teichm{\"u}ller modular groups $\mm$ of Borel-Serre manifolds $X$ 
is well known since the work of Teichm{\"u}ller. 
It is nothing else but the Teichm{\"u}ller spaces $\ttg$\ffdot
Teichm{\"u}ller modular group $\mm$ acts on $\ttg$ discontinuously,
and a subgroup of finite index acts freely by the results of Serre \cite{serre1}.

Motivated by Borel-Serre theory, W. Harvey constructed in \cite{ha2} an analogue of manifolds $\xxb$.
Namely, Harvey 
constructed topological manifolds $\ttb$ (with boundary)
such that\vspace*{-\smallskipamount} 
\[
\ttg\eeq\ttb\aasetmiss\partial\ttb\dff.
\]

\vspace*{-\bigskipamount}\vspace*{-\smallskipamount}
In other words, $\ttg$ is the interior of $\ttb$\nsp.\trs
Both $\ttb$ and $\partial\ttb$ are non-compact.\dss
The boundary $\partial\ttb$ is called the\sss \emph{Harvey boundary} of $\ttg$\nsp.\trs
The canonical action of $\mm$ on $\ttg$ extends to $\ttb$ by the continuity.
This extended action has the following properties:
\vspace{-\bigskipamount}
\begin{itemize}
\item[\phantom{ii}(i)] \hspace*{2.4em}\emph{the action is properly discontinuous}\hff;
\item[\phantom{i}(ii)] \hspace*{2.4em}\emph{a subgroup of finite index in\dss $\mm$\dss acts on\dss $\ttb$\dss freely}\hff;
\item[(iii)] \hspace*{2.4em}\emph{the quotient space\dss $\ttb\dff/\mm$\dss is compact}\hff. 
\end{itemize}

\vspace{-0.2ex}
\vspace{-\bigskipamount}
In particular, $\ttb\dff/\mm$ is a compactification of the moduli space $\ttg/\mm$\dfdot
In fact, $\ttb$ is not only a topological manifold;\dss it has a structure of a smooth manifold with corners.
This structure is not completely canonical (a subtle choice is involved in its construction;\dss see \cite{i-harvey}). 
But any natural construction of such a structure leads to a $\mm$\nsp-invariant structure.
Therefore, we may assume that it is $\mm$\nsp-invariant.
Then for any subgroup $\Gamma$ of $\mm$ acting freely on $\ttb$ 
the quotient $\ttb\dff/\Gamma$ is a smooth manifold with corners.

We refer to the paper \cite{bs} by A. Borel and J.-P. Serre for the definition and the basic properties of manifolds with corners. 
Since the theory of the Harvey boundary is to a big extent modeled on the theory of A. Borel and J.-P. Serre \cite{bs},
this seems to be the most natural reference. 
In the present paper we will need only one result of the theory of manifolds with corners;
see the proof of the next Lemma.

\mypar{Lemma.}{free-quotients} \emph{Suppose that\dss $\Gamma$ is a subgroup of\trs $\mm$\dss of finite index in\sss $\mm$\nsp.\trs
If\dss $\Gamma$\sss acts freely on\sss $\ttb$, then\sss $\ttb\dff/\fff\Gamma$\sss is finitely triangulable space of type\sss $K(\Gamma, 1)$\dnsp.}

\proof Since every topological manifold with boundary is homotopy equivalent to its interior,
$\ttb\trf/\dff\Gamma$\sss is homotopy equivalent to $\ttg\fff/\fff\Gamma$\dfdot
Since $\ttg\fff/\fff\Gamma$ is a $K(\Gamma, 1)$\dnsp-space, as we mentioned in Section \ref{introduction},
$\ttb\trf/\dff\Gamma$\sss is also a $K(\Gamma, 1)$\dnsp-space.
In addition, $\ttb\trf/\dff\Gamma$\sss is a smooth manifold with corners.
It is known that the corners of a smooth manifolds with corners can be smoothed (see \cite{bs}).
Therefore, $\ttb\trf/\dff\Gamma$\sss is homeomorphic to a smooth manifold (without corners).
Since $\Gamma$ is a subgroup of finite index in $\mm$,
and $\ttb\trf/\dff\mm$\sss is compact, the quotient $\ttb\trf/\dff\Gamma$\sss is also compact.
As is well known, every compact smooth manifold is finitely triangulable.
Therefore, $\ttb\trf/\dff\Gamma$\sss is finitely triangulable. 
This completes the proof of the lemma.  \eproof

By the property (ii) of the Harvey boundary there is a subgroup $\Gamma$ of $\mm$ which acts on $\ttb$ freely.
By Lemma \ref{free-quotients}, the quotient \sss $\ttb\dff/\fff\Gamma$\sss admits a finite triangulation.
In particular, it is homotopy equivalent to a finite $CW$-complex.
Therefore, the action of $\Gamma$ on $\ttb$ fits into the framework of Theorems \ref{bie-duality} and \ref{vcd-conn-general} (with $\xxb=\ttb$).

\mypar{Lemma.}{vcd-conn} \emph{If the reduced homology groups\dss $H_i(\partial\ttb)\eeq 0$\sss 
for\sss $i\eeq 0\dff,\dff 1\dff,\dff\ldots\dff,\dff c{\minus}1$\dnsp,
then}
\[
\vcd(\mm\hspace*{0.04em})\leq \dim\ttg{\minus}1{\minus}c = 6\genus{\minus}7{\minus}c\endss.
\]

\vspace*{-\medskipamount}
This is a special case of Theorem  \ref{vcd-conn-general}.
The case $c=2$ is sufficient for the applications in this paper,
and we will prove Lemma \ref{vcd-conn} only in this case. 
See Lemma \ref{h0h1-est} below. 
In fact, this proof works mutatis mutandis in the general situation of Theorem \ref{vcd-conn-general}.

Recall that a group $\Gamma$ is called a \emph{group of type} (FL) 
if the trivial $\Gamma$\dnsp-module $\zz$ admits 
a resolution of finite length consisting of finitely generated free $\Gamma$\dnsp-modules.

\mypar{Lemma.}{type-fl} \emph{Under the assumptions of Lemma \ref{free-quotients},
$\Gamma$ is a group of type (FL).}

\proof The lemma follows from Lemma \ref{free-quotients} together with Proposition 9 of \cite{serre}.  \eproof

\mypar{Lemma.}{beieck} \emph{Let $k$ be a natural number. 
If $\Gamma$ is a group of finite cohomology dimension,
and if $H^n(\Gamma,M)\eeq 0$ for and $n\gres k$ and all free $\Gamma$\dnsp-modules $M$\dfcom
then $\ccd\Gamma\leqs k$\ffdot}

\proof This lemma is due to R. Bieri and B. Eckmann \cite{bie}. See \cite{bie}, Proposition 2.1.  \eproof

\mypar{Lemma.}{h0h1-est} \emph{If the reduced homology groups $H_0(\partial\ttb)\eeq H_1(\partial\ttb)\eeq 0$,
then} 
\[
\vcd(\mm)\leqs 6\genus{\minus}9\endss.
\]

\vspace*{-\bigskipamount}
\proof Let $\Gamma$ be a subgroup of finite index of $\mm$.
We may assume that the action of $\Gamma$ on $\ttg$ is free.
It is sufficient to prove that under such assumptions the cohomology dimension $\ccd(\Gamma)\leqs 6\genus{\minus}9$\dfdot
We start with the following claim.

{\sc Claim 1.} \emph{$H^n(\Gamma, \zz\hff[\Gamma]\hff)\eeq 0$ for $n\gres 6\genus{\minus}9$\dfcom
where $\zz\hff[\Gamma]$ is the integer group ring of $\Gamma$ together with its standard structure of a right $\Gamma$\dnsp-module
(\hff given by the multiplication in $\Gamma$\nsp).}

\emph{Proof of the claim.}\hspace*{0.2em}\ Since $\ttb\trf/\dff\Gamma$\sss is finitely triangulable $K(\Gamma,1)$\dnsp-space and
$\ttb$ is its universal covering (because $\ttb$ is homotopy equivalent to $\ttg$
and hence is a contractible space),
we can apply the results of R. Bieri and B. Eckmann \cite{bie} (see \cite{bie}, Subsection 6.4).
Their results imply that
\[
H^n(\Gamma,\zz\hff[\Gamma]\hff) = H_{d{\minus}n{\minus}1} (\partial\ttb\dff,\zz)
\]
for all $n$\dfcom where $d\eeq\dim\ttg$\nsp.\sss
By the assumptions of the lemma, $H_k(\partial\ttb\dff,\zz)\eeq 0$ for $k\less 2$\dfdot
Therefore $H^n(\Gamma,\zz\hff[\Gamma]\hff)\eeq 0$ for $d{\minus}n{\minus}1\less 2$\dfcom
i.e. for $n\gres d{\minus}3$\dfdot
Hence $H^n(\Gamma,\zz\hff[\Gamma]\hff)\eeq 0$ for 
$n\gres d{\minus}3\eeq \dim\ttg{\minus}3\eeq 6\genus{\minus}6{\minus}3\eeq 6\genus{\minus}9$\dfdot
This proves Claim 1.  \esubproof

Recall that the module $\zz\hff[\Gamma]$ is a free $\Gamma$\dnsp-module with one free generator.
The next step is to extend the above claim to arbitrary free modules.

{\sc Claim 2.} \emph{In $M$ is a free $\Gamma$\dnsp-module, then
$H^n(\Gamma, M)\eeq 0$ for $n\gres 6\genus{\minus}9$\dfdot}

\emph{Proof of the claim.}\hspace*{0.2em}\ Since any finitely generated free $\Gamma$\dnsp-module is
isomorphic to a finite sum of copies of $\zz\hff[\Gamma]$\dfcom
Claim 1 implies that $H^n(\Gamma, M)\eeq 0$ for $n\gres 6\genus{\minus}9$ for any finitely generated free module $M$\dfdot
By Corollary \ref{type-fl} the group $\Gamma$ is a group of type (FL).
Therefore, the functors $H^n(\Gamma,\bullet)$ commute with direct limits by \cite{serre}, Proposition 4.
It follows that $H^n(\Gamma, M)\eeq 0$ for $n\gres 6\genus{\minus}9$ and every free $\Gamma$\dnsp-module $M$\dfdot
This proves Claim 2.  \esubproof

It remains to note that the cohomology dimension of $\Gamma$ is finite (see Section \ref{introduction})
and apply Lemma \ref{beieck}. This completes the proof of Lemma \ref{h0h1-est}.  \eproof

\myitpar{Deduction of Theorem \ref{vcd} from  Theorem \ref{ccsimply}.} 
Now we can prove that Theorem \ref{ccsimply} implies Theorem \ref{vcd}.
By a result of W. Harvey \cite{ha2}, the boundary $\partial\ttb\neq\varnothing$\dfdot
By another result of W. Harvey \cite{ha2}, $\partial\ttb$ is homotopy equivalent to $\ccg$\dfdot
By combining Lemma \ref{h0h1-est} with these results of W. Harvey, 
we see that Theorem \ref{ccsimply} implies the first part of Theorem \ref{vcd},
namely, that $\vcd(\mm)\leqs 6\genus{\minus}9$ if $\genus\geqs 2$\dfdot

It remains to prove that Theorem \ref{ccsimply} implies the part of Theorem \ref{vcd} concerned with $\Mod_2$\dfdot
First, note that Theorem \ref{ccsimply} together with Lemma \ref{h0h1-est} imply that 
\begin{equation}
\label{est-above}
\vcd(\Mod_2)\dff\dff\leq\dff\dff {6\cdot2}\dff\dff -\dff\dff 6\dff\dff =\dff\dff 3\endss.
\end{equation}
On the other hand, by Section \ref{introduction}
\begin{equation}
\label{est-below}
3\genus\dff\dff -\dff\dff 3\dff\dff\leq\dff\dff \vcd(\mm)\endss.
\end{equation} 
for all $\genus$\nsp. 
By applying (\ref{est-below}) to $\genus\eeq 2$\dnsp,\dss
we see that $3\leqs \vcd(\Mod_2)$\dfdot
By taking (\ref{est-above}) into account, we see that $\vcd(\Mod_2)\eeq 3$\dfdot
It remains to prove that $\Mod_2$ is virtually a duality group.

\mypar{Theorem}{genus2} \emph{$\Mod_2$ is virtually a duality group.}

\proof Since $\dim \ccc(X_2)\leqs 2$ and $\ccc(X_2)$ is simply-connected by Theorem \ref{ccsimply},
$\ccc(X_2)$ is homotopy equivalent to a wedge of $2$\dnsp-spheres.
Hence $\partial\overline{\tei}\hspace{0.5pt}_2$ is also homotopy equivalent to a wedge of $2$\dnsp-spheres.
In particular, $H_i(\partial\overline{\tei}\hspace{0.5pt}_2)=0$ if $i\nneq 2$\ffcom
and $H_2(\partial\overline{\tei}\hspace{0.5pt}_2)$ is torsion free.
Let $\Gamma$ be a subgroup $\Mod_2$ acting freely on $\overline{\tei}\hspace{0.5pt}_2$ 
and having finite index in $\Mod_2$\dfdot
Replacing, if necessary, $\Gamma$ by a subgroup of index $2$ in $\Gamma$\dfcom
we can assume that the manifold $\overline{\tei}\hspace{0.5pt}_2/\Gamma$ is orientable.
It remains to apply Theorem \ref{bie-duality} to $\Gamma$ and $\xxb=\overline{\tei}\hspace{0.5pt}_2$\dnsp.\dss \eproof

\mysection{The complex of curves and the Hatcher-Thurston complex}{ccht}

\myitpar{Simplicial complexes.} By a \emph{simplicial complex} 
$\mathcal V$ we understand a simplicial complex in the sense of E. Spanier \cite{sp},
i.e. a pair consisting of a set $\mbox{V}$ together with a collection of \emph{finite subsets of}\/ $\mbox{V}$\dfdot 
As usual, we think of $\mathcal V$ as a structure on the set $\mbox{V}$\dfcom
namely a structure of a simplicial complex.
Elements of $\mbox{V}$ are called\dss \emph{the vertices of\/ $\mathcal V$}\dfcom and subsets of $\mbox{V}$ 
from the given collection are called \emph{the simplices of\/ $\mathcal V$}\dfdot
These data are required to satisfy only one condition:\dss \emph{a subset of a simplex is also a simplex}.
The \emph{dimension} of simplex $S$ is defined as $\dim S\eeq(\card S)\hff{\minus}\hff1$\dfdot
The \emph{dimension} of simplicial complex $\mathcal V$ is defined as the maximum of dimensions of its simplices,
if such maximum exists, and as the infinity $\infty$ otherwise.

\myitpar{Geometric realizations.} Every simplicial complex $\mathcal V$ canonically defines a topological space, which is called 
the \emph{geometric realization} of $\mathcal V$ and denoted by $|\fff\mathcal V\dff |$\dfdot
The idea is to take a copy $\Delta_S$ of the standard geometric simplex $\Delta_{\dim S}$ 
for every simplex $S$ of $\mathcal V$\dfcom and to glue simplices $\Delta_S$ together in such a way 
that $\Delta_T$ will be a face of $\Delta_S$ if $T\ssub S$\dfcom
i.e. if $T$ is a \emph{face} of $S$ in the sense of theory of simplicial complexes.
We omit the details.

When we speak about topological properties of simplicial complexes, they should be understood
as properties of the geometric realization.
The main properties of interest for us, namely, the connectedness and the simply-connectedness,
can be defined pure combinatorially in terms of simplicial complexes, but such an approach is
cumbersome and hides the main ideas.

\myitpar{Barycentric subdivisions.} Every simplicial complex $\mathcal V$ canonically defines another simplicial complex, which is
called the \emph{barycentric subdivision} of $\mathcal V$ and denoted by $\mathcal V'$\dfdot
The vertices of $\mathcal V'$ are the simplices of $\mathcal V$\dfdot
A set of vertices of $\mathcal V'$\dfcom
i.e. a set of simplices of $\mathcal V$\dfcom
is a simplex of $\mathcal V'$ if and only if it has the form
$\{\fff S_1\dff,\dff S_2\dff,\dff\ldots\dff,\dff S_n\fff\}$
for some chain $S_1\ssub S_2\ssub\ldots\ssub S_n$ of simplices of $\mathcal V$\dfdot
A vertex $v$ of $\mathcal V$ is usually identified with the $0$\dnsp-dimensional simplex
$\{\hff v\hff\}$ of $\mathcal V$\dfcom and, hence, with a vertex of ${\mathcal V}'$\dfdot

As is well known, taking the barycentric subdivision does not changes the geometric realization.  
In other terms, for any simplicial complex $\mathcal V$ there is a canonical homeomorphism between 
$|\fff{\mathcal V}'\dff |$ and  $|\fff\mathcal V\dff |$\dfdot

\myitpar{Circles on surfaces and their isotopy classes.} As usual, 
we call by a \emph{simple closed curve} on a surface $X$ (not necessarily closed) a one-dimensional closed connected submanifold of $X$\dfdot
A simple closed curve on a surface $X$ is also called a \emph{circle} on $X$\dfdot
For a circle $C$ in a surface $X$ we will denote \emph{the isotopy class of\/ $C$ in\/ $X$} by $\llvv{C}$\dfdot
The surface $X$ is usually clear from the context, even if $C$ is also a circle in some
other relevant surfaces (for example, some subsurfaces of $X$).
For a collection $C_1\dff,\dff C_2\dff,\dff\ldots\dff,\dff C_n$ we will denote by $\llvv{C_1\dff,\dff C_2\dff,\dff\ldots\dff,\dff C_n}$
the set of the isotopy classes $\llvv{C_i}$ of circles $C_i$ with $1\leqs i\leqs n$\dfdot
\newcommand{\lrlr}[1]{{\langles{\hspace*{0.05em}{#1}\hspace*{0.05em}}\rangles}}
In other terms,
\[
\langles\hspace*{0.02em} C_1\dff,\dff C_2\dff,\dff\ldots\dff,\dff C_n\hspace*{0.05em}\rangles 
= {\{\hspace*{0.05em}}\lrlr{C_1}\dff,\dff \lrlr{C_2}\dff,\dff\ldots\dff,\dff \lrlr{C_n\hff}{\hspace*{0.05em}\}} 
\]
Recall that a circle on $X$ is called \emph{non-trivial} if it cannot be deformed in $X$ into a point or into a boundary component of $X$\dfdot

\myitpar{Complexes of curves.} If $X$ is a compact surface, possibly with non-empty boundary,
then \emph{complex of curves}\/ $\ccc(X)$ is a simplicial complex in the above sense.
The vertices of $\ccc(X)$ are the isotopy classes $\lrlr{C}$ of non-trivial circles $C$ in $X$\dfdot 
A collection of such isotopy classes is a simplex if and only if it is either empty, or
the isotopy classes from this collection can be represented by pair-wise disjoint circles. 
In other words, if $C_1\dff,\dff C_2\dff,\dff\ldots\dff,\dff C_n$ are pair-wise disjoint circles on $X$\dfcom
then the set $\llvv{C_1\dff,\dff C_2\dff,\dff\ldots\dff,\dff C_n}$ is a simplex of $\ccc(X)$\dfcom 
and there are no other simplices (the empty set is the only simplex with $n{\eeq}0$\dnsp;\dss its dimension is $n{\minus}1\eeq 0{\minus}1\eeq {\minus}1$\nsp).

It is well known that if $X$ is a closed orientable surface of genus $\genus$ and\hspace*{0.2em} if 
$C_1\dff,\qff C_2\dff,\qff\ldots\dff,\qff C_n$ are pair-wise disjoint and 
pair-wise non-isotopic circles on $X$\hspace*{-0.2em}, 
then $n\leqs 3\genus{\minus}3$\hspace*{-0.2em},
and there are such collections with $n\eeq 3\genus{\minus}3$ (for example, the collection of circles on\dss Fig. 1).
It follows that $\dim\ccc(X)\eeq 3\genus{\minus}4$ if $X$ is a closed orientable surface of genus $\genus$\ffdot

\myitpar{Hatcher-Thurston complexes.} As before, we denote by $X_{\hff\genus}$ a closed orientable surface of genus $\genus$.
A set $\{\fff C_1\dff,\dff C_2\dff,\dff\ldots\dff,\dff C_{\fff\genus}\trf\}$ of $\genus$ circles on $X_{\hff\genus}$ is called
a \emph{geometric cut system} on $X_{\hff\genus}$ if the circles $C_1\dff,\dff C_2\dff,\dff\ldots\dff,\dff C_{\fff\genus}$ are pair-wise disjoint
and the complement $X_{\hff\genus}\aasetmiss(C_1\cup\ldots\cup C_{\fff\genus})$ is (homeomorphic to) a $2\genus$\dnsp-punctured sphere.
If $\{\fff C_1\dff,\dff C_2\dff,\dff\ldots\dff,\dff C_{\fff\genus}\trf\}$ is a geometric cut system,
then we call the set of the isotopy classes $\llvv{C_1\dff,\dff C_2\dff,\dff\ldots\dff,\dff C_{\fff\genus}}$ a \emph{cut system}.

Suppose that $\{\fff C_1\dff,\dff C_2\dff,\dff\ldots\dff,\dff C_{\fff\genus}\trf\}$ is a geometric cut system on $X_{\hff\genus}$\dfdot
Suppose that $1\leqs i\leqs\genus$, and that $C'$ be a circle on $X_{\hff\genus}$ disjoint from circles $C_j$ with $j\nneq i$\dfcom
and transversely intersecting $C_i$ at exactly $1$ point.
If we replace $C_i$ by $C'$ in $\{\fff C_1\dff,\dff C_2\dff,\dff\ldots\dff,\dff C_{\fff\genus}\trf\}$\dfcom
we get another cut system. A \emph{simple move} is the operation of replacing the geometric cut system 
\[
\{\fff C_1\dff,\dff\ldots\dff,C_i\dff,\dff\ldots\dff,\dff C_{\fff\genus}\trf\}
\hspace*{0.3em}\mbox{ by the geometric cut system }\hspace*{0.3em}
\{\fff C_1\dff,\dff\ldots\dff,C'\dff,\dff\ldots\dff,\dff C_{\fff\genus}\trf\}\endss,
\]
and also the corresponding operation of replacing the cut system
\[
\llvv{C_1\dff,\dff\ldots\dff,C_i\dff,\dff\ldots\dff,\dff C_{\fff\genus}}
\hspace*{0.3em}\mbox{ by the cut system }\hspace*{0.3em}
\{\llvv{C_1\dff,\dff\ldots\dff,C'\dff,\dff\ldots\dff,\dff C_{\fff\genus}}\nsp.
\]
{\em Usually we will describe a simple move by pictures omitting the unchanging circles.}

Some sequences of simple moves are \emph{cycles} in the sense that they begin and end at the same geometric cut system.
The three special types of cycles, depicted on Fig. 3, are the key ingredients of the construction of the Hatcher-Thurston complexes.
It is assumed that the circles omitted from the pictures are disjoint from the ones presented, and form cut systems with them. 

The Hatcher-Thurston complex $\cht$ of $X_{\hff\genus}$ is a $2$\dnsp-dimensional cell complex 
(it is not a simplicial complex) constructed as follows.
Every cut system 
$\llvv{C_1\dff,\dff C_2\dff,\dff\ldots\dff,\dff C_{\fff\genus}}$
is a $0$\dnsp-cell of $\cht$\dnsp;\dss
there are no other $0$\dnsp-cells.
If one $0$\dnsp-cell can be obtained from another by a simple move, then these two $0$\dnsp-cells
are connected by a $1$\dnsp-cell corresponding to this move;
there are no other $1$\dnsp-cells.
Clearly, one geometric cut system can be obtained from another one by no more than one simple move,
and even a cut system can be obtained from another one by no more than one simple move.
Therefore, two $0$\dnsp-cells of $\cht$ are connected by no more than one $1$\dnsp-cell.
At this moment we have already a $1$\dnsp-dimensional cell complex consisting of the just
described $0$\dnsp-cells and $1$\dnsp-cells.
It is denoted by $\chtone$\dfdot 
The Hatcher-Thurston complex $\cht$ is obtained from $\chtone$ by attaching $2$\dnsp-cells to $\chtone$
along circles resulting from the three special types of cycles, namely, cycles of types (I), (II), and (III).

This definition of $\cht$ was suggested by A. Hatcher and W. Thurston \cite{ht}, 
who also proved the following fundamental result.

\renewcommand{\topfraction}{1.0}
\renewcommand{\textfraction}{0.0}

\begin{figure}[p]
\includegraphics[width=0.98\textwidth]{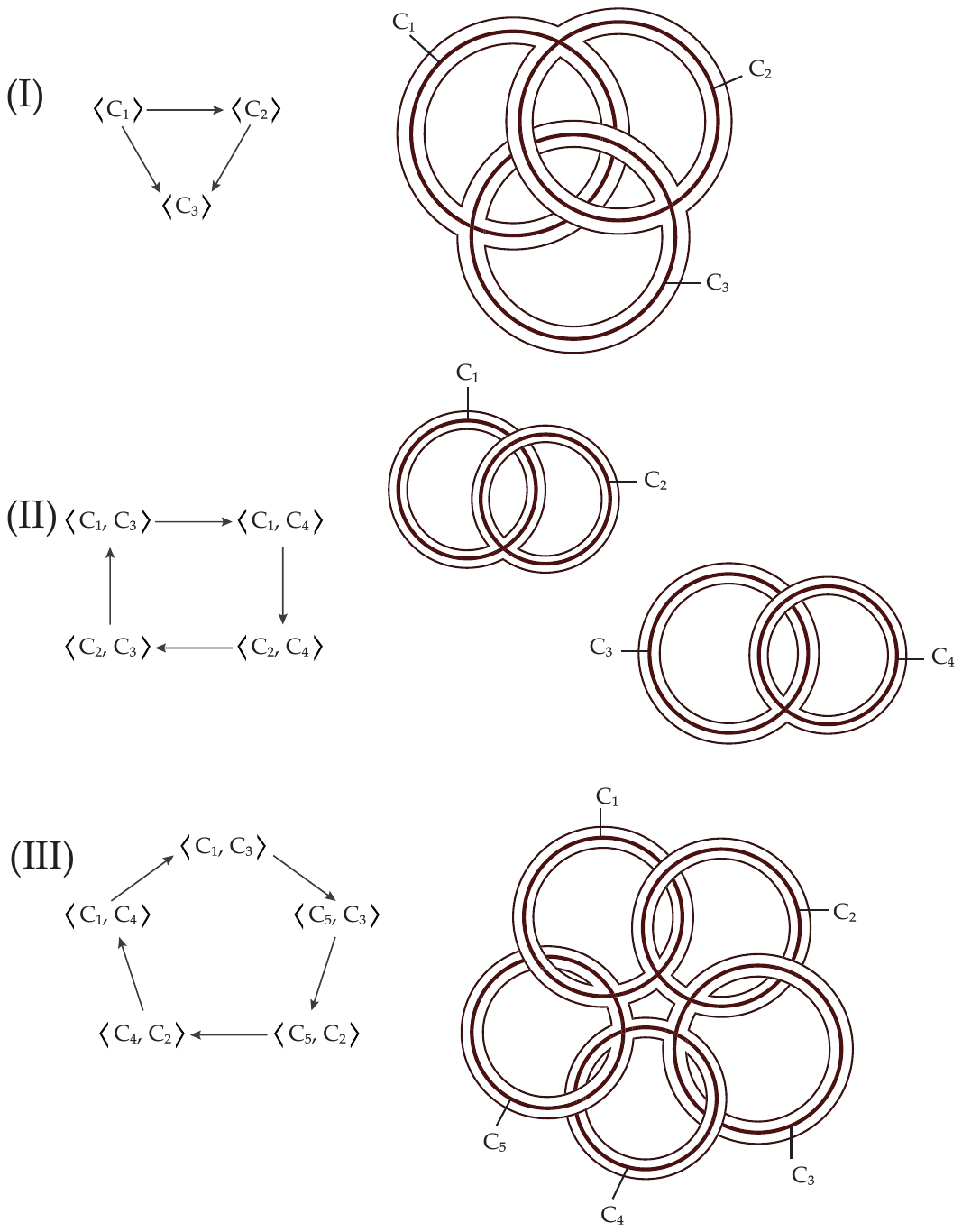}
\caption{The Hatcher-Thurston moves.}
\label{moves}
\end{figure}

\mypar{Theorem.}{ht-theorem} \emph{If the genus $\genus$ of $X_{\hff\genus}$ is $\nsp\geqs 2$\ffcom\sss
then the cell complex $\cht$ is connected and simply-connected.}

\myitpar{Proof of Theorem \ref{ccsimply}.} It is based on a construction of
a map ${\mathcal J}\colon\chtone\tto | \ccc(\xx)\trf |$ such that the following two lemmas hold.

\mypar{Lemma.}{extension} \emph{If\dss $\genus\geqs 2$\ffcom\sss
then map ${\mathcal J}$ can be extended to a map $\cht\tto|\ccc(\xx)\trf|$\dfdot}

\mypar{Lemma.}{loops} \emph{If\dss $\genus\geqs 2$\ffcom
then every loop in $|\ccc(\xx)\trf |$ is freely homotopic to a loop of the form ${\mathcal J}\hff(\beta)$\dfcom
where $\beta$ is a loop in $\chtone$\dfdot}

Since $\cht$ is simply connected by Theorem \ref{ht-theorem}, Lemmas \ref{extension} and \ref{loops} together
imply that $\ccc(\xx)$ is simply connected.
Therefore, the proof of Theorem \ref{ccsimply} is completed modulo Lemmas \ref{extension} and \ref{loops} 
and the construction of ${\mathcal J}$\ffdot\vspace*{\medskipamount}

{\em In the rest of this section we assume that $\genus\geqs 2$\ffdot}

\myitpar{The construction of ${\mathcal J}$\ffdot} The cell complex $\chtone$ is the geometric realization
of a simplicial complex ${\mathcal R}\hspace{0.5pt}(\xx)$ defined as follows.
The set of vertices of ${\mathcal R}\hspace{0.5pt}(\xx)$ is equal to the set of $0$\dnsp-cells of $\cht$\dfdot
In other words, the vertices of ${\mathcal R}\hspace{0.5pt}(\xx)$ are the cut systems on $\xx$.
If two $0$\dnsp-cells $\mbox{V}_1$\ffcom $\mbox{V}_2$ of $\cht$ are connected by a $1$\dnsp-cell of $\cht$ 
(i.e. if they are related by a simple move),
then the pair $\{\hff\mbox{V}_1,\mbox{V}_2\hff\}$ is a simplex of ${\mathcal R}\hspace{0.5pt}(\xx)$ (of dimension $1$\nsp). 
There are no other simplices; in particular, there are no simplices of dimension $\nsp\geqs 2$\dfdot

Let ${\mathcal R}\eeq{\mathcal R}\hspace{0.5pt}(\xx)$\dfcom $\ccc\eeq\ccc(\xx)$\dfdot
Let ${\mathcal R}'$\dfcom $\ccc'$ be the barycentric subdivisions of ${\mathcal R}$\dfcom $\ccc$ respectively.
Let us construct a morphism of simplicial complexes $J\colon{\mathcal R}'\tto\ccc'$\dfdot
Recall that a \emph{morphism of a simplicial complexes} ${\mathcal A}\tto{\mathcal B}$ is defined as a map
of the set of vertices of $\mathcal A$ to the set of vertices of $\mathcal B$ such that the image of a simplex is also a simplex.
If $Z\eeq\llvv{\dnsp C_1\dff,\dff C_2\dff,\dff\ldots\dff,\dff C_{\fff\genus}\fff}$ is a vertex of ${\mathcal R}$
considered as a vertex of ${\mathcal R}'$\dfcom
we set
\[
J\hff (Z) = \llvv{C_1\dff,\dff C_2\dff,\dff\ldots\dff,\dff C_{\fff\genus}\fff}\endss,
\] 
where the right hand side is a simplex of $\ccc$ considered as a vertex of $\ccc'$\dfdot
If $Z\eeq\{\hff \mbox{V}_1,\mbox{V}_2\hff\}$ is the vertex of ${\mathcal R}'$ 
corresponding to to the edge of $\mathcal R$ connecting the vertices $\mbox{V}_1$ and $\mbox{V}_2$ of ${\mathcal R}$\dfcom
then these two vertices connected by a simple move.
Let
\[
\mbox{V}_1 = \llvv{C_1\dff,\dff\ldots\dff,\dff C_i\dff,\dff\ldots\dff C_{\fff\genus}\fff}
\mapsto
\llvv{C_1\dff,\dff\ldots\dff,\dff C_i'\dff,\dff\ldots\dff C_{\fff\genus}\fff} =\fff \mbox{V}_2
\]
be this simple move.
Then we set 
\[
J\hff (Z) = \llvv{C_1\dff,\dff\ldots\dff,C_{i{\minus}1}\dff,\dff C_{i{\plus}1}\dff,\dff\ldots\dff,\dff C_{\fff\genus}}
\]
where, again, the right hand side is a simplex of $\ccc$ considered as a vertex of $\ccc'$\dfdot
Obviously,\dss the map $J$ from the set of vertices of ${\mathcal R}'$ to the set of vertices of $\ccc'$ is
a morphism of simplicial complexes ${\mathcal R}'\tto\ccc'$\dfdot 

Recall that a morphism of simplicial complexes $f\colon{\mathcal A}\tto{\mathcal B}$ 
canonically defines a continuous map $|\dff f\dff |\colon|\fff{\mathcal A}\fff|\tto|\fff{\mathcal B}\fff|$\nsp,\dss
called the \emph{geometric realization} of $f$\dfdot 
Therefore, $J$ leads to a continuous map $|\fff J\fff |\colon |\dff {\mathcal R}'\dff |\tto|\fff \ccc'\fff |$\dfdot
Since the geometric realization of the barycentric subdivision of a complex is canonically homeomorphic to the
geometric realization of the complex itself, we may consider $|\fff J\fff |$ as a map 
$|\dff {\mathcal R}\dff |\tto|\fff \ccc\fff |$\dfdot

Recall that $\chtone$ is the geometric realization of ${\mathcal R}\eeq{\mathcal R}\hspace{0.5pt}(\xx)$\dfdot
Therefore, we may define ${\mathcal J}\colon\chtone\tto | \ccc(\xx)\trf |$ to be the map $|\dff J\dff |$ considered as
map $\chtone\tto | \ccc(\xx)\trf |$\dfdot

\mypar{Proof of Lemma \ref{extension}\hff.}{proof-ext} It is sufficient to prove that ${\mathcal J}$ maps every cycle of type (I), (II), or (III)
to a loop contractible in ${\mathcal C}$\dfdot
We will concider these three types of cycles separately.

\myitpar{Cycles of type \textnormal{\textup{(I)}}.} Let $\llvv{C_1\dff,\dff\ldots\dff,\dff C_{\fff\genus}}$
be a $0$\dnsp-cell of $\cht$ involved into a cycle of type (I).
Then all circles $C_1\dff,\dff\ldots\dff,\dff C_{\fff\genus}$ except one remain unchanged under $3$ simple moves 
forming this cycle.
We may assume that the circle  $C_1$ is not changing.
Then $J$ maps every vertex of ${\mathcal R}'$ which belongs to the geometric realization of this cycle into a vertex of ${\mathcal C}'$ 
of the form $\llvv{C_1\dff,\dff\ldots\dff\ldots}$\dfdot
Therefore $J$ maps every such vertex into a vertex of ${\mathcal C}'$ contained in the star of the vertex $\llvv{C_1}$ of ${\mathcal C}$\dfcom
considered as a vertex of ${\mathcal C}'$ (by a standard abuse of notations, we identify $\llvv{C_1}$ with $\{\llvv{C_1}\}$).\dss 
Indeed, $\{ \llvv{C_1}, \llvv{C_1\dff,\dff\ldots\dff\ldots\dff}\}$ is an edge of ${\mathcal C}'$ connecting
$\{ \llvv{C_1} \}$ with $\{ \llvv{C_1}, \llvv{C_1\dff,\dff\ldots\dff\ldots\dff}\}$\dfdot
It follows that the image of this cycle under the morphism $J$ is contained in the star of $\llvv{C_1}$
and hence the image of the circle resulting from this cycle under ${\mathcal J}$ is contained in the geometric realization of this star,
and hence is contractible in this geometric realization.
Therefore it is contractible in $|\fff \ccc'\fff |\eeq|\fff \ccc\fff |\eeq|\fff \ccc(\xx)\trf |$\dfdot 
This completes the proof for the cycles of type (I).  \esubproof

\myitpar{Cycles of type \textnormal{\textup{(II)}}.} Let us consider a cycle of type (II).
The simple moves of such a cycle change two circles, 
and the other $\genus{\minus}2$ circles do not change.
Therefore, if $\genus\geqs 3$\dfcom then at least one circle of the cut systems from this cycle remains in place under all $4$ simple moves of this cycle.
This allows to complete the proof in this case in exactly the same way as we dealt with the cycles of type (I).

It remains to consider the case of $\genus\eeq 2$\dfdot
In this case each cut system consists of $2$ circles and only four circles $C_1\dff,\dff C_2\dff,\dff C_3\dff,\dff C_4$ 
are involved in the cycle.\dss 
See Fig. \ref{moves} (II).\dss
In this case there exist a non-trivial circle $C_0$ on $X$ disjoint from $C_1\dff,\dff C_2\dff,\dff C_3\dff,\dff C_4$\nsp.\sss
For example, the union $C_1\cup C_2$ is contained in a subsurface of $\xx$ diffeomorphic to a torus with one hole.
We can take as $C_0$ the boundary circle of this torus with one hole.
Alternatively, we can define $C_0$ as the circle dividing $\xx$ into two tori with $1$ boundary component each
such that $C_1\cup C_2$ is contained in one of them, and $C_3\cup C_4$ is contained in the other one.
(This more symmetric description of $C_0$ easily implies that $C_0$ is unique up to isotopy, but we will not need this fact.)

Every $0$\dnsp-cell of $\mathcal{HT}(\xx)$ occurring in our cycle  
has the form $\llvv{C_i\dff,\dff C_j\fff}$\dfcom
where $i\eeq 1\mbox{ or }2$ and $j\eeq 3\mbox{ or }4$\dfdot
In order to describe this cycle in more details, it is convenient to introduce an involution $\sigma$
on the set $\{\hff 1\dff,\dff 2\dff,\dff 3\dff,\dff 4\fff \}$\dfdot
Namely, we set
\[
\sigma(1) = 2\fff,\quad \sigma(2) = 1\fff,\quad \sigma(3) = 4\fff,\quad \sigma(4) = 3\fff.
\]
Then every $1$\dnsp-cell contained in our cycle corresponds to a simple move of the form
\[
\llvv{C_i\dff,\dff C_j\fff}\mapsto \llvv{C_{\sigma(i)}\dff,\dff C_j},\hspace*{0.5em}
\mbox{ or of the form }\hspace*{0.5em}\llvv{C_i\dff,\dff C_j\fff}\mapsto \llvv{C_i\dff,\dff C_{\sigma(j)}},
\]
where $i\eeq 1\mbox{ or }2$\dfcom and $j\eeq 3\mbox{ or }4$\dfdot
In the barycentric subdivision ${\mathcal R}'$ the edge connecting $\llvv{C_i\dff,\dff C_j}$ with $\llvv{C_{\sigma(i)}\dff,\dff C_j}$
is subdivided into two edges, connecting the vertex 
\[
\{\fff\llvv{C_i\dff,\dff C_j},\dff \llvv{C_{\sigma(i)}\dff,\dff C_j}\fff\}
\]
of ${\mathcal R}'$ with the vertices $\llvv{C_i\dff,\dff C_j}$ and $\llvv{C_{\sigma(i)}\dff,\dff C_j}$ respectively.
Since $C_0$ is disjoint from the circles $C_1\dff,\dff C_2\dff,\dff C_3\dff,\dff C_4$\nsp,\sss
the images of both these edges under the map $J$ are contained in the star of the vertex $\llvv{C_0}$ 
(more precisely, $\{\hspace*{0.025em}\llvv{C_0}\hspace*{0.025em}\}$) of ${\mathcal C}'$\dfdot
The same argument applies to all edges into which our cycle is subdivided in ${\mathcal R}'$\dfdot
It follows that $J$ maps the subdivided cycle into the star of $\llvv{C_0}$ in ${\mathcal C}'$\dfcom
and hence the geometric realization ${\mathcal J}\eeq |\dff J\dff |$ maps the geometric realization of our cycle
into the geometric realization of this star.
Therefore, this image is contractible in the geometric realization of this star, and hence
in $|\fff \ccc'\fff |\eeq|\fff \ccc\fff |\eeq|\fff \ccc(\xx)\trf |$\nsp. 
This completes the proof for the cycles of type (II).  \esubproof

\myitpar{Cycles of type \textnormal{\textup{(III)}}.} This is the most difficult case.
If $\genus\geqs 3$\dfcom then one of the circles is not changed under all five moves of the cycle and 
we can use the same argument as we used for the cycles of type (I) and for the cycles of type (II) in the case $\genus\geqs 3$\dfdot

\begin{figure}[ht]
\includegraphics[width=0.96\textwidth]{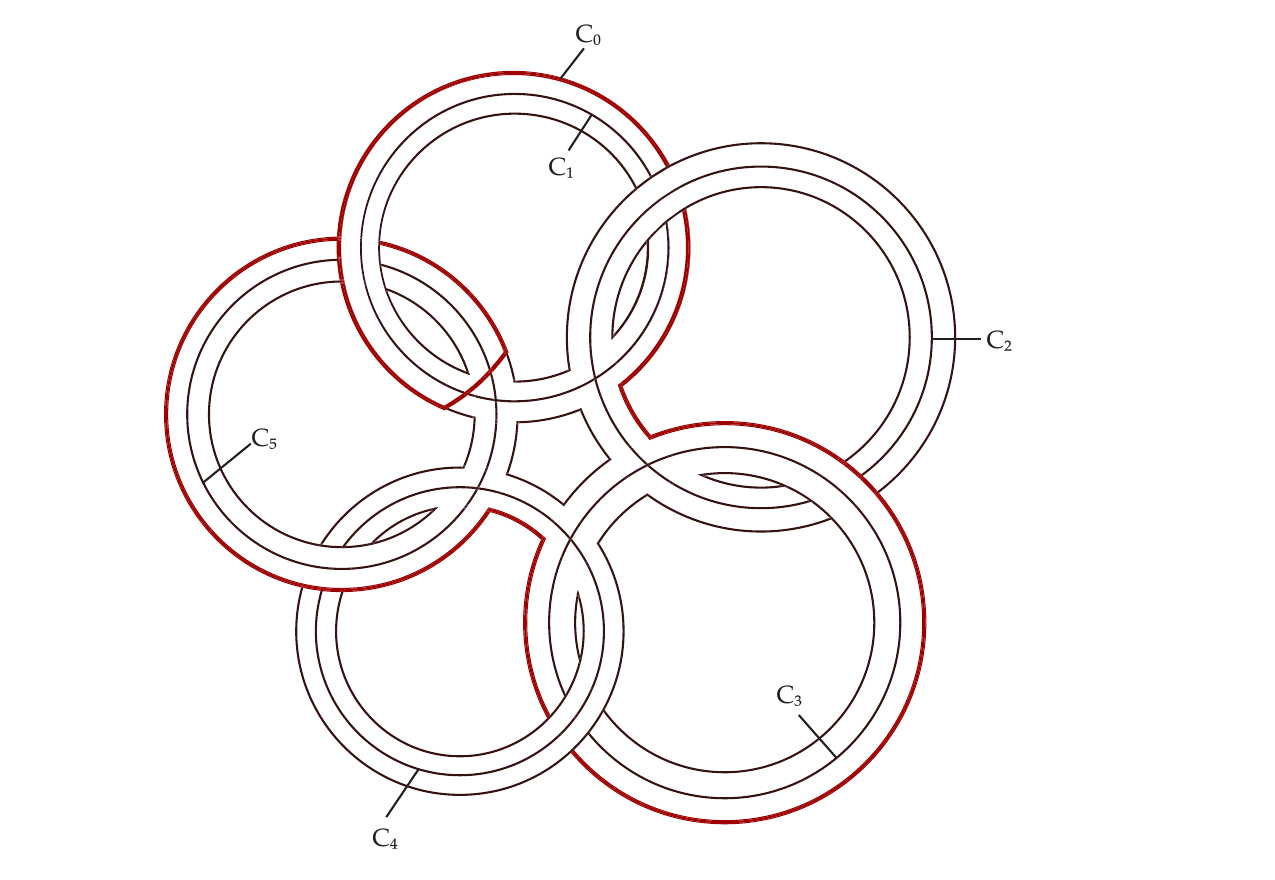}
\caption{Auxiliary circle for move (III).}
\label{3}
\end{figure}

It remains to consider the case of $\genus\eeq 2$\dfdot
In this case each cut system consists of $2$ circles and only five circles $C_1\dff,\dff C_2\dff,\dff C_3\dff,\dff C_4\dff,\dff C_5$ 
are involved in the cycle.\dss 
See Fig. \ref{moves} (III).\dss
Let $C_0$ be a circle on $\xx$ disjoint from $C_2\dff,\dff C_3\dff,\dff C_4$ and intersection each of the circles $C_1$ and $C_5$
transversely at one point. 
One can take as $C_0$ the circle $C_0$ on the Fig. \ref{3}.

\renewcommand{\topfraction}{1.0}
\renewcommand{\textfraction}{0.0}

\renewcommand{\topfraction}{1.0}
\renewcommand{\textfraction}{0.0}

\begin{figure}[ht]
\includegraphics[width=0.96\textwidth]{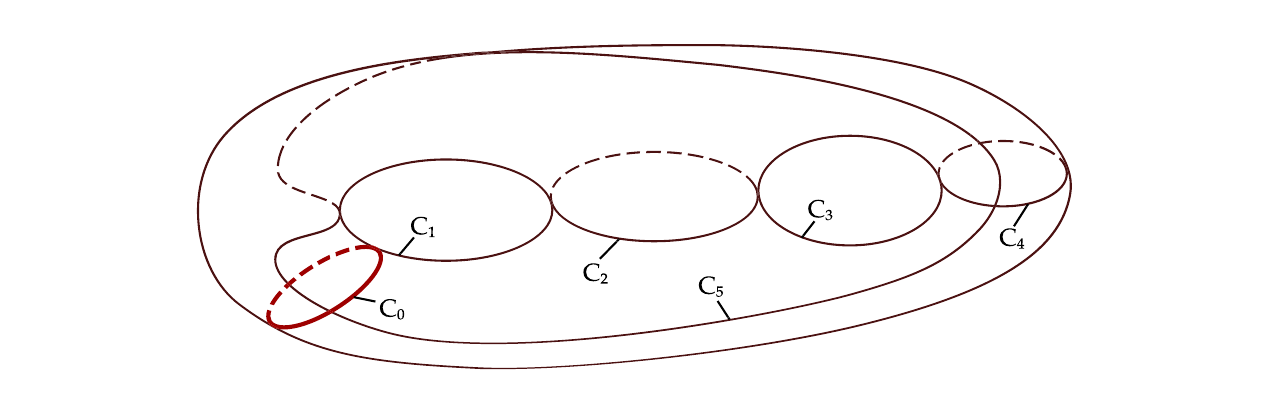}
\caption{Another view of the auxiliary circle for move (III).}
\label{4}
\end{figure}

An alternative way to draw such a circle is presented on the  Fig. \ref{4}.
We leave to the interested readers to show that the circles $C_0$ on these two pictures are isotopic\fff;\sss
we will not use this fact.

\begin{figure}[h]
\includegraphics[width=\textwidth]{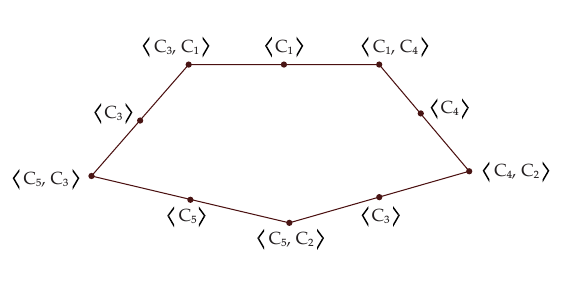}
\caption{The pentagon.}
\label{pentagon}
\end{figure}

Let us consider the image under $\mathcal J$ of the circle in $\chtone$ resulting from our cycle.
This image is the geometric realization of the (simplicial) loop in ${\mathcal C}'$ shown on Fig. \ref{pentagon}.

\renewcommand{\topfraction}{1.0}
\renewcommand{\textfraction}{0.0}

\begin{figure}[h]
\includegraphics[width=\textwidth]{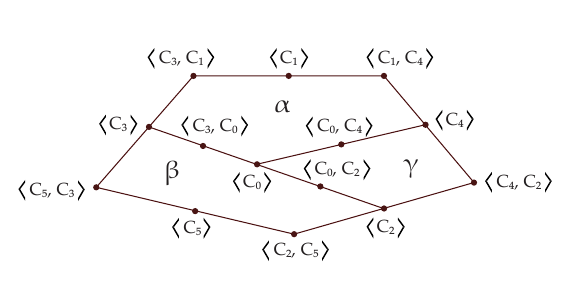}
\caption{Filling in the pentagon.}
\label{filled-pentagon}
\end{figure}

The subgraph (i.e. a $1$\dnsp-dimensional simplicial subcomplex) of ${\mathcal C}'$\dfcom
shown on Fig. \ref{filled-pentagon}
contains the above simplicial loop as a subgraph.

\renewcommand{\topfraction}{1.0}
\renewcommand{\textfraction}{0.0}

\begin{figure}[h]
\includegraphics[width=\textwidth]{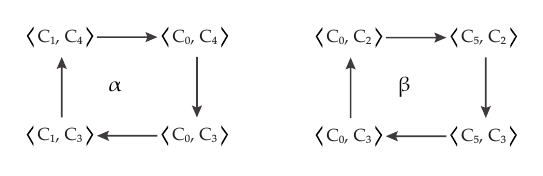}
\caption{Cycles for $\alpha$ and $\beta$\dfdot}
\label{alpha-beta}
\end{figure}

The subgraphs bounding the domains $\alpha$ and $\beta$ on this picture are equal 
to the images under the map $J$ of the barycentric subdivisions
of the two cycles of type (II) in ${\mathcal R}$ shown on Fig. \ref{alpha-beta}.
Therefore, their geometric realizations are contractible in $|\ccc\fff|\eeq |\ccc(\xx)\trf|$\dfdot

\renewcommand{\topfraction}{1.0}
\renewcommand{\textfraction}{0.0}

The subgraph bounding the domain $\gamma$ on this picture is equal 
to the image under the map $J$ of the barycentric subdivision
of the boundary of the triangle (i.e. a $2$\dnsp-dimensional simplex) 
$\llvv{C_0\dff,\dff C_2\dff,\dff C_4\hff}$ in ${\mathcal C}$\dfdot
Therefore, its geometric realization is contractible in $|\ccc\fff|$\dfdot

It follows that the geometric realizations of these $3$ loops (subgraphs) are contractible in $\cht$\dfdot
Therefore, the geometric realization of the loop on Fig. 5 is also contractible in $|\ccc\fff|$\dfdot
Since this geometric realization is the image of our cycle of type (III),\qss
this completes the proof for cycles of type (III), and hence the proof of the lemma.  \esubproof  \eproof

\mypar{Proof of Lemma \ref{loops}\hff.}{proof-loops} Every loop in $| \fff {\mathcal C}\fff | \eeq | \fff \ccc(\xx)\fff |$
is freely homotopic to the geometric realization of a simplicial loop in the $1$\dnsp-skeleton of $\mathcal C$\dfdot
A simplicial loop in the $1$\dnsp-skeleton of $\mathcal C$ is just a sequence
\begin{equation}
\label{generic-loop}
\hspace*{0em}\llvv{C_1},\dff\hspace*{0.3em} \llvv{\dnsp C_2},\dff\hspace*{0.3em} \ldots\ldots\dff,\dff\hspace*{0.3em} \llvv{\dnsp C_n}
\end{equation}
of vertices of $\mathcal C$ such that $\llvv{C_i}$ is connected by an edge of $\mathcal C$ with $\llvv{C_{i{\plus}1}}$
for all $i\eeq 1\dff,\dff 2\dff,\dff\ldots\dff,\dff n{\minus}1$ and 
$\llvv{C_n}$ is connected by an edge with $\llvv{C_1}$\dfdot

\emph{From now on we will interpret $n{\plus}1$ as $1$\dfdot}

Without loss of generality we may assume that $\vvv{\qff C_i\nsp}\nsp \nneq \vvv{\dff C_{i{\plus}1} }$ for
every $i\eeq 1\dff,\dff 2\dff,\dff\ldots\dff,\dff n$\dfdot

{\sc Claim 1.} \emph{Without loss of generality,\dss we can assume that circles $C_i$ are non-separating.}

\emph{Proof of the claim.} Suppose that $C_i$ is a separating circle. Let $Y_0$ and $Y_1$ be two subsurfaces of $\xx$ into
which $C_i$ divides $\xx$\dfdot
Since $\xx$ is a closed surface, both $Y_0$ and $Y_1$ are surfaces with one boundary component resulting from $C_i$\dfdot
Since $C_i$ is a non-trivial circle, neither $Y_0$\dfcom nor $Y_1$ is a disc.
Hence each of surfaces $Y_0$ and $Y_1$ has genus $\geqs 2$\dfdot

Let us first consider the case when both circles $C_{i{\minus}1}$ and $C_{i{\plus}1}$ are non-separating.

If $C_{i{\minus}1}$ and $C_{i{\plus}1}$ are contained in the same part $Y_j$ of $\xx$ (where $j\eeq 0\mbox{ or }1$\nsp),
then we can choose a non-separating circle $C'_i$ in the other part $Y_{1{\minus}j}$ of $\xx$,
because $Y_{1{\minus}j}$ is a surface of genus $\geqs 2$\dfdot
Then both
\[
\llvv{C_{i{\minus}1}\dff,\dff C_i\dff,\dff C_{i{\plus}1}}\hspace*{0.5em}\mbox{ and }
\hspace*{0.5em}\llvv{C_{i{\minus}1}\dff,\dff C'_i\dff,\dff C_{i{\plus}1}}
\]
are simplices (triangles). 
Therefore, our loop is homotopic to the loop resulting from replacing $\llvv{C_i}$ by $\llvv{C'_i}$ in it.
Since the circle $C'_i$ is non-separating in $Y_{1{\minus}j}$\dfcom 
it is non-separating in $\xx$,
and hence our new simplicial loop involves one separating circle less than the original one.

If $C_{i{\minus}1}$ and $C_{i{\plus}1}$ are contained in different parts of $\xx$,
then $C_{i{\minus}1}\cap C_{i{\plus}1}\eeq\varnothing$\dfcom 
and hence $\llvv{ C_{i{\minus}1} }$ and $\llvv{ C_{i{\plus}1} }$ are connected by an edge in $\mathcal C$\dfdot
Moreover, $\llvv{ C_{i{\minus}1}\dff,\dff C_i\dff,\dff C_{i{\plus}1} }$
is a simplex (triangle) of $\mathcal C$\dfdot
It follows that if we delete $\llvv{C_i}$ from our loop, we get a new loop
which is homotopic to the original one.
As before, the new loop involves one separating circle less than the original one.

Let us now consider the case when the circle $C_{i{\plus}1}$ is separating (and $C_i$ is also separating, as before).
We may assume that the circles $C_i$ and $C_{i{\plus}1}$ are disjoint (replacing them by isotopic circles, if necessary).
Then $C_i$ and $C_{i{\plus}1}$ together divide $\xx$ into three parts 
$Z_0\dff,\dff Z_1\dff,\dff Z_2$\dfdot
Since the circles $C_i$ and $C_{i{\plus}1}$ are non-isotopic (by our assumption) and are both non-trivial,
each of the surfaces $Z_0\dff,\dff Z_1\dff,\dff Z_2$ has genus $\geqs 1$\dfdot
Since the circle $C_{i{\minus}1}$ is disjoint from $C_i$\dfcom
the circle $C_{i{\minus}1}$ may intersect no more than two of surfaces $Z_0\dff,\dff Z_1\dff,\dff Z_2$\dfdot
Let $Z_k$ be a part disjoint from $C_{i{\minus}1}$\dfdot 
Let $C'_i$ be some non-separating circle in $Z_k$ (such a circle exists because the genus of $Z_k$ is $\geqs 1$\nsp).
Then the circles $C_{i{\minus}1}\dff,\trf C'_i\dff,\trf C_{i{\plus}1}$ are pair-wise disjoint,
and hence both
\[
\llvv{C_{i{\minus}1}\dff,\dff C_i\dff,\dff C_{i{\plus}1}}\hspace*{0.5em}\mbox{ and }
\hspace*{0.5em}\llvv{C_{i{\minus}1}\dff,\dff C'_i\dff,\dff C_{i{\plus}1}}
\]
are simplices (triangles).
It follows that if we replace in our loop the vertex $\llvv{C_i}$ by the vertex $\llvv{C'_i}$\dfcom
we will get a new loop homotopic to the original one.
Since $C'_i$ is non-separating circle in $Z_k$\dfcom
and hence is a non-separating circle in $\xx$,
the new loop involves one separating circle less than the original one.

Finally, in the case when the circle $C_{i{\minus}1}$ is separating, the same arguments as in the case when
$C_{i{\plus}1}$ apply.
This allows us replace our loop by a homotopic new loop involving one separating circle less than the original one in this case also.

By repeating the above procedure until there will be no separating circles involved, we can
construct a new loop homotopic to the original loop and involving no separating circle.
This proves our claim.  \esubproof

{\sc Claim 2.} \emph{Without loss of generality,\dss we can assume that, in addition to circles $C_i$ being non-separating,
every edge $\vv{\hff C_i\dff,\dff C_{i{\plus}1}\dnsp}\eeq\fff\{\dff\vvv{\hff C_i},\dff \vvv{C_{i{\plus}1}\dnsp}\hff\}$ can be completed to a cut system.}

\emph{Proof of the claim.} By Claim 1, we can assume that all circles $C_i$ are non-separating.
Suppose that $\llvv{C_i\dff,\dff C_{i{\plus}1}\nsp}$ cannot be completed to a cut system.
Since $\llvv{C_i}$ and $\llvv{C_{i{\plus}1}\nsp}$ are connected by an edge of $\mathcal C$\dfcom
we may assume that the circle $C_i$ and $C_{i{\plus}1}$ are disjoint.
Then $\llvv{C_i\dff,\dff C_{i{\plus}1}\nsp}$ cannot be completed to a cut system only if
the union $C_i\cup C_{i{\plus}1}$ divides our surface $\xx$ into two parts 
(it cannot divide $\xx$ into three parts because neither $C_i$\dfcom nor $C_{i{\plus}1}$ divide $\xx$).
Let these two parts be $Y_0$ and $Y_1$, so that $Y_0\cup Y_1\eeq\hff\xx$ and
$Y_0\cap Y_1\eeq\hff C_i\cup C_{i{\plus}1}$\dfdot
Since $\llvv{C_i}$ and $\llvv{C_{i{\plus}1}}$ are assumed to be different, and hence
$C_i$ and $C_{i{\plus}1}$ are not isotopic, each of the subsurfaces $Y_0$ and $Y_1$ has genus $\geqs 1$\dfdot

Let us choose some circle $C'_i$ contained in $Y_0$ and non-separating in $Y_0$
(this is possible because the genus of $Y_0$ is $\geqs 1$\nsp).
Then $C'_i$ is non-separating in $\xx$ also and the circles $C_i\dff,\dff C'_i\dff,\dff C_{i{\plus}1}$ are pair-wise disjoint.
In particular,  $\llvv{C_i\dff,\dff C'_i\dff,\dff C_{i{\plus}1}\nsp}$ is $2$\dnsp-simplex (triangle) of $\mathcal C$\dfdot
Moreover, since $C'_i$ is non-separating in $Y_0$\nsp,\dss it is also non-separating in both $\xx\setmiss C_i$ and $\xx\setmiss C_{i{\plus}1}$\dfdot
Therefore both unions $C_i\cup C'_i$ and $C'_i\cup C_{i{\plus}1}$ do not divide $\xx$ into two parts. 
It follows that both pairs $\llvv{C_i\dff,\dff C'_i}$ and $\llvv{C'_i\dff,\dff C_{i{\plus}1}\nsp}$ can be completed to cut systems.

Let us replace the edge connecting $\llvv{C_i}$ with $\llvv{C_{i{\plus}1}\nsp}$ in our loop
by the following two edges: the first one connecting $\llvv{C_i}$ with $\llvv{C'_i}$\dnsp;\dss
the second one connecting $\llvv{C'_i}$ with $\llvv{C_{i{\plus}1}\nsp}$.
Since $\llvv{C_i\dff,\dff C'_i\dff,\dff C_{i{\plus}1}}$ is $2$\dnsp-simplex, the new loop is homotopic to the original one.
Since both pairs $\llvv{C_i\cup C'_i}$ and $\llvv{C'_i\cup C_{i{\plus}1}\nsp}$ can be completed to cut systems, 
the new loop has less edges which cannot be completed to cut systems than the original one.

By repeating this procedure we can construct a new loop homotopic to the original one and having the required properties.
This completes the proof of the  claim.  \esubproof

By Claim 2, it is sufficient to consider loops (\ref{generic-loop}) in $\mathcal C$ such that every $C_i$ is a non-separating circle,
and every pair $\llvv{C_i\dff,\dff C_{i{\plus}1}\nsp}$ can be extended to a cut system.
Given such a loop (\ref{generic-loop}),
we consider the following loop in the barycentric subdivision ${\mathcal C}'$
\begin{equation}
\label{sub-loop}
\hspace*{0em}\llvv{C_1},\dff\hspace*{0.3em}\llvv{C_1\dff,\dff C_2},\dff\hspace*{0.3em} \llvv{\dnsp C_2},\dff\hspace*{0.3em}\llvv{C_2\dff,\dff C_3},\dff\hspace*{0.3em} \ldots\ldots\dff,\dff\hspace*{0.3em}\llvv{C_{n{\minus}1}\dff,\dff C_n},\dff\hspace*{0.3em} \llvv{\dnsp C_n}\dff,\dff\hspace*{0.3em}\llvv{C_n\dff,\dff C_1}
\end{equation}
For every $i\eeq 1\dff,\dff 2\dff,\dff\ldots\dff,\dff n$ there is an edge of this loop connecting $\llvv{C_i}$ with $\llvv{C_i\dff,\dff C_{i{\plus}1}\nsp}$\dnsp,\dss
and an edge connecting $\llvv{C_i\dff,\dff C_{i{\plus}1}\nsp}$ with $\llvv{C_{i{\plus}1}\nsp}$ (recall that $n{\plus}1$ is interpreted as $1$\nsp);
there are no other edges.
Clearly, the loops (\ref{generic-loop}) and (\ref{sub-loop}) have the same geometric realization. 

Let us complete each pair $\llvv{C_i\dff,\dff C_{i{\plus}1}\nsp}$ to a cut system $Z_i$\dfdot
Clearly, $Z_i$ has the form $Z_i\dff\eeq\llvv{\dnsp C_i\dff,\dff C_{i{\plus}1}\dff,\dff C^i_3\dff,\dff\ldots\dff,\dff C^i_{\fff\genus}\fff}$
if $\genus\geqs 3$\dfcom and $Z_i\dff\eeq\llvv{\dnsp C_i\dff,\dff C_{i{\plus}1}\nsp}$ if $\genus\eeq 2$\dfdot

Let us temporarily fix an integer $i$ between $1$ and $n$\dfdot
Let us cut our surface $\xx$ along $C_i$ and denote the result by $X^0_i$\dfdot
The surface $X^0_i$ has two boundary components and its genus is equal to $\genus{\minus}1$\dfdot
Next, let us glue two discs to the components of $\partial X^0_i$ and denote the result by $X^1_i$\dfdot
Clearly, $X^1_i$ is a closed surface of genus $\genus{\minus}1$\dfdot
We may assume that $\{\fff C_i\dff,\dff C_{i{\plus}1}\dff,\dff C^i_3\dff,\dff\ldots\dff,\dff C^i_{\fff\genus}\trf\}$
and $\{\fff C_{i{\minus}1}\dff,\dff C_{i}\dff,\dff C^{i{\minus}1}_3\dff,\dff\ldots\dff,\dff C^{i{\minus}1}_{\fff\genus}\trf\}$
are geometric cut systems on $\xx$. 
Then  $\{\fff  C_{i{\plus}1}\dff,\dff C^i_3\dff,\dff\ldots\dff,\dff C^i_{\fff\genus}\trf\}$
and $\{\fff C_{i{\minus}1}\dff,\dff C^{i{\minus}1}_3\dff,\dff\ldots\dff,\dff C^{i{\minus}1}_{\fff\genus}\trf\}$
are geometric cut systems on $X^1_i$ (because all circles involved are contained in $X^0_i\ssub X^1_i$\nsp).
Therefore, by taking the isotopy classes in $X^1_i$ instead of $\xx$,
we can define two cut system on $X^1_i$ as follows:
$Z^0_i\dff\eeq\llvv{\dnsp C_{i{\plus}1}\dff,\dff C^i_3\dff,\dff\ldots\dff,\dff C^i_{\fff\genus}}$
and $Z^1_i\dff\eeq\llvv{\dnsp C_{i{\minus}1}\dff,\dff C^{i{\minus}1}_3\dff,\dff\ldots\dff,\dff C^{i{\minus}1}_{\fff\genus}}$\dfdot
Because $\mathcal{HT}\hff(X^1_i)$ is connected,
$Z^0_i$ can be joined with $Z^1_i$ by a path in $\mathcal{HT}\hff(X^1_i)$\dfcom
and hence by a path in the $1$\dnsp-skeleton of $\mathcal{HT}\hff(X^1_i)$\dfdot
It follows that $Z^0_i$ can be joined with $Z^1_i$ by a path in the simplicial complex ${\mathcal R}'(X^1_i)$\dfdot
Let us denote this path by $\alpha_i$\dfdot
Let $\llvv{D_1\dff,\dff D_2\dff,\dff\ldots\dff,\dff D_n}$ be a vertex of $J(\alpha_i)$
(where $n\eeq\genus{\minus}1$ or $\genus{\minus}2$ by the construction of $J$\nsp)\hff.\sss
Because $X^1_i\dff\setmiss\inter X^0_i$ is a union of two disjoint discs\dfcom
we may assume, replacing the circles $D_1\dff,\dff D_2\dff,\dff\ldots\dff,\dff D_n$ by circles isotopic to them in $X^1_i$\dfcom
if necessary, that $D_1\dff,\dff D_2\dff,\dff\ldots\dff,\dff D_n\ssub\inter X^0_i$.\sss
Then $\llvv{C_i\dff,\dff D_1\dff,\dff D_2\dff,\dff\ldots\dff,\dff D_n}$ is a vertex of $\ccc'(\xx\hff)$\dfdot
By adding in this way $\llvv{C_i}$ to all vertices of the path $\alpha_i$\dfcom
we will obtain a sequence of vertices of $\ccc'\eeq{\mathcal C}'(\xx\hff)$\dfdot
Clearly, this sequence is a simplicial path in $\ccc'$\dfcom
and, moreover, it is equal to $J(\beta_i)$ for some simplicial path $\beta_i$ in ${\mathcal R}'(\xx\hff)$\dfdot

Now, let us put together all paths $\beta_i$ for $i\eeq 1\dff,\dff 2\dff,\dff\ldots\dff,\dff n$\dfdot
Let $\beta$ be the resulting loop.
In order to complete the proof, it is sufficient to show that the geometric realization of the loop (\ref{sub-loop})
is freely homotopic to $J(\beta)$\dfdot
In order to prove this, it is sufficient, in turn, to prove that for every $i$ the path $J(\beta_i)$\dfcom
which connects $Z_{i{\minus}1}$ with $Z_{i}$\dfcom
is homotopic relatively to the endpoints to the path
\[
\llvv{Z_{i{\minus}1}},\dff\hspace*{0.3em}\llvv{C_{i{\minus}1}\dff,\dff C_i},\dff\hspace*{0.3em}\llvv{C_i},\dff\hspace*{0.3em}
\llvv{C_i\dff,\dff C_{i{\minus}1}},\dff\hspace*{0.3em}\llvv{Z_i}\endss.
\]
But both these paths are contained in the star of $\llvv{C_i}$ in $\ccc'$\dfdot
Therefore they are homotopic relatively to the endpoints.
This completes the proof of the lemma.  \eproof

Lemmas \ref{loops} and \ref{extension} are now proved.
As we saw, these lemmas imply Theorem \ref{ccsimply}.
In addition, Theorem \ref{vcd} follows from Theorem \ref{ccsimply} by the results of Section \ref{boundary}.
Therefore, our main theorems, namely, Theorems \ref{vcd} and \ref{ccsimply} are proved.

\mysection{Beyond the simply-connectedness of $\ccc(\xx\hspace*{0.04em})$}{context}

\vspace*{\bigskipamount}
\myitpar{The connectedness and simply-connectedness of $\ccc(\xx\hspace*{0.04em})$\dnsp.} The connectedness 
of $\ccc(\xx\hspace*{0.04em})$ can be proved by a direct argument, 
which we leave as an exercise to an interested reader.

A natural approach to proving the simply-connectedness of $\ccc(\xx\hspace*{0.04em})$ 
is to look for a reduction of this problem to the simply-connectedness of $\cht$\dnsp.\qss 
The latter was proved by A. Hatcher and W. Thurston \cite{ht}.\dss  
The complexes $\ccc(\xx\hspace*{0.04em})$ and $\cht$ are not related in any direct and obvious manner.
Still, it is possible to relate them in a not quite direct (but canonical) way and 
deduce the simply-connectedness and connectedness of $\ccc(\xx\hspace*{0.04em})$ from the corresponding properties of $\cht$\dfdot
This deduction is the heart of the paper \cite{iv-vcd} and is presented in Section \ref{ccht} above.
This deduction allows to prove that $\ccc(\xx\hspace*{0.04em})$ is simply-connected if $\genus\geqs 2$
(note that $\ccc(\xx\hspace*{0.04em})$ is not even connected if $\genus\eeq 1$\nsp). 

Suppose that $\genus\geqs 2$\dfdot
In view of the results of Sections \ref{introduction} and \ref{boundary}, 
the connectedness of $\ccc(\xx\hspace*{0.04em})$ 
implies an estimate of $\vcd(\mm\hspace*{0.04em})$ better than 
the trivial estimate $\vcd(\mm\hspace*{0.04em})\leqs 6\genus{\minus}7$\dnsp,\qss
and the simply-connectedness implies an even better estimate.
Namely, the connectedness of $\ccc(\xx\hspace*{0.04em})$ implies that $\vcd(\mm\hspace*{0.04em})\leqs 6\genus{\minus}8$\dfcom
and the simply-connectedness implies that $\vcd(\mm\hspace*{0.04em})\leqs 6\genus{\minus}9$\dfdot

\myitpar{The complexes of curves and the Hatcher-Thurston complexes.} At the first sight, 
deducing the simply-connectedness of
$\ccc(\xx)$ from the simply-connectedness of $\cht$ seems to be somewhat artificial.
This was my opinion in 1983 and for many years to follow.
Much later, with the benefit of the hindsight, I started to think that this opinion was short-sighted.
In fact, this deduction contains the nuclei of many arguments used later to study the complexes of curves, 
starting with my papers \cite{iv5}, \cite{i-stab}, and \cite{i-ihes,i-conf,i-imrn}.
Nowadays these arguments are among the most natural tools of trade.

There was also a better reason to be unsatisfied with such a deduction. 
Namely, such a deduction cannot be extended to prove higher connectivity of $\ccc(\xx\hspace*{0.04em})$ when it is expected,
since the complex $\cht$ is only $2$\dnsp-dimensional.
The idea to generalize the whole paper of A. Hatcher and W. Thurston \cite{ht} 
to a higher-dimensional complexes, yet to be constructed, appeared to be too far-fetched.
This is the \emph{road not taken}.

A natural alternative to constructing higher-dimensional versions of $\cht$\dfcom
is to try to apply the ideas of \cite{ht} to $\ccc(\xx\hspace*{0.04em})$ directly.
The main tool of Hatcher and Thurston \cite{ht} is the Morse-Cerf theory \cite{cerf},
an analogue of the Morse theory for families of functions with $1$ parameter.
In order to work with the complex of curves one needs, first of all, 
to modify the Morse-Cerf theory in such a way that it will lead to result about
$\ccc(\xx\hspace*{0.04em})$\dfcom and not about $\cht$\dfdot
Also, one needs to at least partially extend the Morse-Cerf theory 
to the families of functions on surfaces with arbitrary number of parameters.
The latter would be necessary even if the high-dimensional versions of $\cht$ would be constructed.

\myitpar{The classification of singularities and the Morse-Cerf theory.} The Morse theory deals with individual functions,
which may be considered as families of functions with $0$ parameters. 
J. Cerf \cite{cerf} extended the Morse theory to families of functions with $1$ parameter.
Families with $2$ parameters also appear in \cite{cerf}, but they are not arbitrary:
they are constructed in order to deform families with $1$ parameter.
The main difficulty in extending the Morse-Cerf theory to families with an arbitrary number of parameters
results from the lack of classification of singularities of functions in generic families of functions depending on several parameters.
The Morse theory requires only the classification of singularities of generic functions.
The Cerf theory \cite{cerf} requires only the classification of singularities of functions in generic families with $1$ parameter.

The term \emph{classification} is used here in a precise and very strong sense.
A \emph{singular point} of a smooth function $f\colon M\tto\rr$\dfcom where $M$ is a smooth manifold, 
is defined as a point $x\elem M$ such that the differential $d_x f$ is equal to $0$\dfdot
Two singular points $x$\dfcom $y$ of functions $f\colon M\tto\rr$\dfcom $g\colon N\tto\rr$ respectively are called \emph{equivalent}\dss
if $g\circ\varphi\eeq f{\pplus}c$ for some real constant $c$ and some diffeomorphism $\varphi$ 
between a neighborhood of $x$ in $M$ and a neighborhood of $y$ in $N$\dfcom
such that $\varphi(x)\eeq y$\dfdot
A \emph{singularity} can be defined as an equivalence class of singular points of smooth functions.
A \emph{classification of singularities} of function in some class consists of a list of all possible singularities of functions in this class,
and explicit formulas for representatives of each singularity in this list.
An explicit formula for a representative is called a \emph{normal form} of the corresponding singularity.
Usually one takes as a normal form of a singularity a polynomial in several variables with $0$ being the singular point in question.

The most important classes of functions for the applications are the classes of functions occurring in \emph{generic} families
of functions with a given number $m$ of parameters.
The singularities of such functions are called the singularities\sss \emph{of codimension} $\leqs m$\dfdot
The \emph{codimension of a singularity} is defined as the smallest $m$ such that the singularity is of codimension $\leqs m$\dfdot
The codimension of a singularity is the main measure of its complexity from the point of view of applications.

By 1983 V. I. Arnold and his students to a big extent completed Arnold's program of classification of
singularities of functions (of course, the nature of Arnold's program is such that it never can be completed).
The book \cite{avg}, presenting the main results of Arnold's program, appeared in 1982.

Arnold discovered, in particular, that the codimension is not the best measure of complexity of a singularity for the purposes of classification.
Instead, the dimension of the space of deformations of a singularity is a more appropriate characteristic. 
This dimension, properly defined, is called the \emph{modality} of a singularity.
The singularities of the modality $0$ are, essentially, the ones which cannot be deformed to a non-equivalent singularity by a small deformation.
They are called \emph{simple singularities}. 
The singularities of modality $1$ are called \emph{unimodal}.  
The simple singularities are simplest to classify, the next case being the unimodal singularities.
Arnold classified the simple singularities already in 1972 \cite{a72a}, \cite{a72b} 
(see also \cite{a75}), and the unimodal ones in 1974 \cite{a73}.

As a corollary of the classification of simple singularities, Arnold found a classification all singularities of codimension $\leqs 5$ (they are all simple).
The classification of unimodal singularities lead to a classification of singularities of codimension $\leqs 9$ (they are all either simple or unimodal).
But there is no hope to find a classification (in the above sense) of singularities of arbitrary codimension (or, what is the same, of arbitrary modality).

As it eventually turned out, in the context of our problem the Morse functions are the worst ones,
and one can bypass the classification of singularities entirely.
This was done in \cite{iv5}.
The next section tells more about the story behind \cite{iv5}.

\mysection{Reminiscences\fff:\sss $\vcd(\mm\hspace*{0.05em})$\sss in Leningrad\hff,\sss 1983}{myremin}

\vspace*{\bigskipamount}
\myitpar{The virtual cohomology dimension\sss $\vcd(\mm\hspace*{0.04em})$ and the connectivity of\dss $\ccc(\xx\hspace*{0.04em})$\dnsp.} 
In the Spring of 1983 I was working, among other things, on the problem of
computing the virtual cohomology dimension of $\mm$\ffdot
The arguments of Sections \ref{introduction} and \ref{boundary} were in my mind from the very beginning,
despite the fact that I learned the Bieri--Eckmann theory \cite{bie} only in the process of working on this problem,
and I was only vaguely familiar with the Borel--Serre theory \cite{bs}.
It seems that these ideas were in the air at the time.

Since I admired the Hatcher--Thurston paper \cite{ht} and studied it in details,
it was only natural to try to deduce the simply-connectedness of $\ccc(\xx\hspace*{0.04em})$
from the main result of \cite{ht}, the simply-connectedness of $\cht$\dnsp.\dss
This lead to the arguments of Section \ref{ccht} and a proof of Theorem \ref{ccsimply}. 
In turn, this immediately lead to a proof of Theorem \ref{vcd}.

{\small
This work was done in March and April of 1983 and presented at Rokhlin's Topology Seminar in Leningrad in April.
In May I prepared a research announcement which included Theorems \ref{vcd} and \ref{ccsimply}, 
as well as other results about Teichm{\"u}ller modular groups which I proved starting from December of 1982.
The announcement was presented by Academician L.D. Faddeev to \emph{Doklady of Academy of Sciences of the USSR} 
(known also as \emph{DAN}\fff) at May 16, 1983.
It was published \cite{iv0} in the first months of the next year.
These results were also included in \emph{Short Communications} distributed
at least among the participants of the Warsaw Congress in August of 1983.

} 

\vspace*{-0.8\medskipamount}
\myitpar{The Morse-Cerf theory and the complexes of curves.} Eventually it turned out that there is a method to
apply an ideal version of the Morse-Cerf theory directly to $\ccc(\xx\hspace*{0.04em})$ 
without using the Hatcher-Thurston complex $\cht$ as an intermediary.
I found such a method in Summer of 1983.
In fact, it turned out that it is much easier to apply the Morse-Cerf theory directly to $\ccc(\xx\hspace*{0.04em})$
than to the Hatcher-Thurston complex $\cht$\dfcom
not to say about using $\cht$ as an intermediary.

As expected, the method allowed in principle to prove that the complex of curves $\ccc(\xx\hspace*{0.04em})$ is $n$\dnsp-connected
if a classification of singularities up to codimension $n+1$ is available.
The method required that the normal forms of these singularities were not too complicated in a precise sense.
The well known normal forms of singularities of codimension $\leqs 2$ are trivially not too complicated in this sense.
This allowed to reprove the connectedness and the simply-connectedness of $\ccc(\xx\hspace*{0.04em})$
without using the Hatcher-Thurston theory.
The relation with the classification of singularities was completely parallel to the Hatcher-Thurston theory\hff:
in order to prove that $\cht$ is connected (respectively, simply-connected), Hatcher and Thurston used
the classification of singularities of codimension $\leqs 1$ (respectively, of codimension $\leqs 2$\nsp). 

Arnold's classification of singularities of codimension $\leqs 5$ immediately implied that these
singularities are simple enough for my method to work.
This allowed to prove that $\ccc(\xx\hspace*{0.04em})$ is $3$\dnsp-connected if 
$\genus\geqs 3$\dfcom and is $4$\dnsp-connected if $\genus\geqs 4$\dfdot
In view of Lemma \ref{vcd-conn}, this implied that $\vcd(\mm)\leqs 6\genus{\minus}11$ if $\genus\geqs 3$\dfcom
and $\vcd(\mm)\leqs 6\genus{\minus}12$ if $\genus\geqs 4$\dfdot
After checking the properties of the normal forms of  singularities of codimension $\leqs 6$\dfcom
I proved that, moreover, $\ccc(\xx\hspace*{0.04em})$ is $5$\dnsp-connected if $\genus\geqs 4$\dfdot
In view of Lemma \ref{vcd-conn}, this implied that $\vcd(\mm)\leqs 6\genus{\minus}13$ if $\genus\geqs 4$\dfdot

I planned to go further through Arnold's lists of normal forms, and, in particular, to look at all singularities of codimension $\leqs 9$\dfdot
It was clear that such a straightforward approach relying on normal forms will exhaust its potential soon. 
But the experience with the normal forms at the initial part of Arnold's list lead me to believe 
that all singularities of high codimension are very simple for the purposes of
my method, and that there should be a way to bypass the normal forms and the classification.
This work was interrupted by a trip to Warsaw to attend the Warsaw Congress.

\myitpar{Warsaw Congress, August 1983.} By the time of the Warsaw Congress
I had proved that $\mathrm{vcd}\hspace*{0.1em}(\mm)\leqslant 6g{\minus}11$ if $\genus\geqs 3$
and $\vcd(\mm)\leqs 6\genus{\minus}13$ if $\genus\geqs 4$\dfdot 
It was clear that the method does not stop there. 
I was thrilled when W. Thurston showed up for my short talk at the Congress. 
Unfortunately, my command of spoken English was negligible, and I spoke in Russian. 
Volodya Turaev acted as an interpreter. 
After the talk Thurston suggested to discuss my talk and to tell me the news related to results and problems discussed in my talk. 
Note that at the time the communication between Western and Soviet mathematicians was anything but easy,
and the Warsaw Congress presented a unique opportunity to learn about thing not yet published or even not written down.
During this discussion (with Volodya Turaev continuing to serve as an interpreter)
Thurston told that J. Harer computed the virtual cohomological dimension of $\mm$\dfdot 
Unfortunately, Thurston forgot the actual value of $\vcd(\mm\hspace*{0.04em})$\dfdot 
After being pressed, Thurston agreed that the value of $\vcd(\mm\hspace*{0.04em})$ is {\em ``as expected''}. 

For me, the value of $\vcd(\mm\hspace*{0.04em})$ being {``as expected''} 
meant that everything is parallel to the Borel-Serre theory \cite{bs}.
In particular, $\ccc(\xx\hspace*{0.04em})$ is homotopy equivalent to a bouquet of spheres of dimension equal 
to the topological dimension of $\ccc(\xx\hspace*{0.04em})$\dfcom i.e. to $3\genus{\minus}4$ 
(cf. Remark 3.5 in \cite{ivj}\fff). 
If this is the case, then the Bieri--Eckmann theory \cite{bie} implies that 
$\vcd(\mm\hspace*{0.04em})\eeq\dim \ttg{\minus}(3\genus{\minus}4){\minus}1$\dfdot
Since $\dim \ttg = 6\genus{\minus}\genus$\hspace*{-0.05em},\sss
this means that $\vcd(\mm\hspace*{0.04em})\eeq(6\genus{\minus}6){\minus}(3\genus{\minus}4){\minus}1\eeq 3\genus{\minus}3$\dfdot 
In fact, the Bieri-Eckmann theory \cite{bie}, together with the simply-connectedness of $\ccc(\xx\hspace*{0.04em})$\dfcom  
implies that $\vcd(\mm\hspace*{0.04em})\eeq 3\genus{\minus}3$ if and only if $\ccc(\xx\hspace*{0.04em})$ is
$(3\genus{\minus}5)$\dnsp-connected, but not $(3\genus{\minus}4)$\dnsp-connected (for $\genus\geqs 2$\nsp).

My methods were clearly not sufficient to prove that $\ccc(\xx\hspace*{0.04em})$ is $(3\genus{\minus}5)$\dnsp-connected 
(which is not surprising, because it is indeed not $(3\genus{\minus}5)$\dnsp-connected), 
and I abandoned the project for a couple of months. 

\myitpar{A misunderstanding.} I am inclined to think that W. Thurston wasn't at fault when he said that the value of $\vcd(\mm\hspace*{0.04em})$ 
is {\em ``as expected''} and did not remembered the correct formula. 
W. Thurston was thinking about deeper issues than a formula for $\vcd(\mm\hspace*{0.04em})$\dnsp.\dss
Most likely, he was thinking about the \emph{reasons allowing to find\/}\sss 
the value of the virtual cohomology dimension of  $\mm$\hspace*{-0.05em},\sss and they were {\em ``as expected''}.
My reasons to expect that $\vcd(\mm\hspace*{0.04em})\eeq 3\genus{\minus}3$ were based on an analogy between 
Teichm{\"u}ller modular groups and arithmetic groups.
As it seems now, I expected this analogy to hold with more details than it actually holds.
Since $3\genus{\minus}3$ is equal to the maximal rank of the abelian (and of the solvable) subgroups of 
$\mm$, 
the analogy with the arithmetic groups suggested that $3\genus{\minus}3$ should be the answer.
While this analogy is a very good guiding principle, it is not complete.
Moreover, this lack of completeness makes the theory of Teichm{\"u}ller modular groups much more interesting
than it would be otherwise.

\myitpar{Autumn of 1983.} After returning from Warsaw, I wrote to J. Birman, asking, in particular, 
about what exactly was proved by J. Harer about the order of connectedness of $\ccc(\xx\hspace*{0.04em})$
and the virtual cohomology dimension $\vcd(\mm\hspace*{0.04em})$\dfdot
At that time crossing the USSR border usually took one-two months for a letter.\dss 
The reply from J. Birman arrived only at the late autumn of 1983.\dss
In her reply she wrote me that according to J. Harer  $\ccc(\xx\hspace*{0.04em})$
is $(2\genus{\minus}3)$\dnsp-connected, but is not $2\genus{\minus}2$\dnsp-connected
and that $\vcd(\mm\hspace*{0.04em})\eeq 4\genus{\minus}5$ for $\genus\geqs 2$\dfdot

I immediately realized that this is exactly what my methods can in principle provide.
Independently of the form the classification of singularities takes in higher codimension,
higher than $(2\genus{\minus}3)$\dnsp-connectivity could not be proved by my methods 
because already Morse functions prevent this.
After this I quickly proved that all singularities of higher codimension are indeed simpler 
than the Morse singularities for the purposes of my method.
See \cite{iv5}, Subsection 2.1 and Lemma 2.2 for the key idea. 
This allowed me to complete the proof of $(2\genus{\minus}3)$\dnsp-connectedness of $\ccc(\xx\hspace*{0.04em})$ by the end of 1983.\sss

About the same time the preprint of \cite{ha} arrived. 
It contained, in particular, a beautiful combinatorial proof of the fact that $\ccc(\xx\hspace*{0.04em})$ is homotopy
equivalent to a $(2\genus{\minus}2)$\dnsp-dimensional CW-complex. 
This result is independent from the main part of \cite{ha}, which is concerned with
$(2\genus{\minus}3)$\dnsp-connectedness of $\ccc(\xx\hspace*{0.04em})$\dfdot
Together with the $(2\genus{\minus}3)$\dnsp-connected\-ness of $\ccc(\xx\hspace*{0.04em})$\dfcom
this result implies that $\vcd(\mm\hspace*{0.04em})\eeq\fff 4\genus{\minus}5$ for $\genus\geqs 2$\dnsp.\dss 
Combined with my proof of the $(2\genus{\minus}3)$\dnsp-connectedness of $\ccc(\xx\hspace*{0.04em})$\dnsp,\dss
this leads to a computation of $\vcd(\mm\hspace*{0.04em})$ largely independent from Harer's one.

Harer's exposition was somewhat obscure for my taste, 
and I found a different version of his proof of homotopy equivalence of $\ccc(\xx\hspace*{0.04em})$
to a $(2\genus{\minus}2)$\dnsp-dimensional CW-complex.
It brings to the light the fact that the basic properties of the Euler characteristic
(never mentioned by Harer) are behind Harer's combinatorial arguments.

All these results and their analogues for non-orientable surfaces were published in \cite{iv5}. 

\myitpar{A lemma in Harer's paper.} Harer's paper \cite{ha} contains at least one gap: 
the proof of Lemma 3.6 is not correct and, I believe, cannot be saved.
But the lemma is true. 

In order to state this lemma, let $X_{\hff\genus\fff,\fff s}$ be a surface of genus $\genus$
with $s$ boundary components $\partial_1 X\dff,\dff \partial_2 X\dff,\dff \ldots\dff,\dff \partial_s X$\dfdot
Let us assume that\sss $X_{\hff\genus\fff,\fff 0}\eeq \xx$ and that
$X_{\hff\genus\fff,\fff s{\minus}1}$ results from $X_{\hff\genus\fff,\fff s}$ by glueing a disc
to $\partial_s X$\dfdot
Then every circle in $X_{\hff\genus\fff,\fff s}$ can be considered as a circle in
$X_{\hff\genus\fff,\fff s{\minus}1}$\dfdot
Some circles non-trivial in $X_{\hff\genus\fff,\fff s}$ will be trivial in
$X_{\hff\genus\fff,\fff s{\minus}1}$\dfdot
Let $\widehat{\ccc}(X_{\hff\genus\fff,\fff s}\hspace*{0.04em})$
be the subcomplex of
$\ccc(X_{\hff\genus\fff,\fff s}\hspace*{0.04em})$
having as its vertices the isotopy classes of circles in $X_{\hff\genus\fff,\fff s}$
remaining non-trivial in
$X_{\hff\genus\fff,\fff s{\minus}1}$\dfcom
with the simplices defined as before.
Harer's Lemma 3.6 claims that\dss \emph{for $s>0$ the obvious map
$\widehat{\ccc}(X_{\hff\genus\fff,\fff s}\hspace*{0.04em})
\tto
\ccc(X_{\hff\genus\fff,\fff s{\minus}1}\hspace*{0.04em})$
is a homotopy equivalence.}

The methods of\qss \cite{iv5}\sss do not depend on this lemma,
and in the first version of this paper I wrote
that I am not aware of any proof of it.
In fact, a proof of this lemma is contained in the proof of\dss
Lemma 2.6\dss from my paper\dss \cite{i-stab}.
Also, the main part of the proof of\dss Lemma 2.5 from\dss \cite{i-stab}\dss is a proof
of a similar statement about some complexes build from arcs connecting boundary circles
to themselves (namely, that, in the notation of\qss \cite{i-stab}, the canonical map
$A'(R)\tto B'(S)$ is a homotopy equivalence).
The proofs are similar and depend on some deep results in the theory of minimal surfaces.

In the first version of this paper I wrote also that it is desirable to find
elementary proofs of\dss Harer's\dss Lemma 3.6 and of\dss Lemma 2.5 from\dss \cite{i-stab}.
After the first version of this paper was posted,
A. Hatcher\dss \cite{hat}\dss wrote to me that an elementary proof of the case 
$s\eeq 0$ of\dss Harer's\dss Lemma 3.6
is contained in his paper\dss \cite{hv}\dss with K. Vogtmann.
See\dss \cite{hv}, Proposition 4.7.
Moreover, the same arguments prove the general case of this lemma,
and, mutatis mutandis, lead to an elementary proof of the corresponding
result about complexes build from arcs and hence of\dss Lemma 2.5 from\dss \cite{i-stab}.


\renewcommand{\refname}{\textnormal{Additional references}}

{}\vspace{-\bigskipamount}

\begin{flushright} 

September 30, 2015

The last subsection\dss --\dss September\phantom{0} 3, 2023 

\vspace{\bigskipamount}

1983:\hspace*{3em}Leningrad Branch of\hspace*{4.7em}\hspace*{2.1em} \\ 
Steklov Mathematical Institute,\hspace*{2.1em}{} \\
Leningrad, USSR\hspace*{6.2em}\hspace*{2.1em}{ } \\
\vspace{\bigskipamount}
Current:\hspace*{1.4em}http:/\!/\hspace*{-0.07em}nikolaivivanov.com \hspace*{3.9em}

E-mail:\hspace*{1.7em} ivanov@msu.edu,\hspace*{8.1em}

nikolai.v.ivanov@icloud.com\hspace*{3.25em}
\end{flushright}

\end{document}